\documentclass[11pt]{article}

\usepackage{amsmath,amssymb,amsfonts, amsthm, amsbsy, mathtools, epsfig}
\usepackage[usenames]{color}
\usepackage{rotating}
\usepackage{enumerate}

\usepackage{verbatim}
\usepackage{caption}
\usepackage{subcaption}
\usepackage{bbm}
\usepackage{xcolor}
\usepackage{bm}
\usepackage{mathrsfs,color,dsfont}

\usepackage{thmtools}
\usepackage{thm-restate}
\usepackage{hyperref}
\usepackage{cleveref}

\newtheorem{thm}{Theorem}

\newtheorem{assumption}{Assumption}

\newtheorem{example}[thm]{Example}

\newtheorem{theorem}{Theorem}
\newtheorem{lemma}[theorem]{Lemma}
\newtheorem{proposition}[theorem]{Proposition}

\newcommand{\RNum}[1]{\uppercase\expandafter{\romannumeral #1\relax}}
\makeatother

\newcommand\tabcaption{\def\@captype{table}\caption}

\definecolor{DSgray}{cmyk}{0,1,0,0}

\newcommand{\E}{\mathbb{E}}

\newcommand{\PP}{\mathbf{P}} 

\begin{document}





\title{Scheduling with Service-Time Information: The Power of Two Priority Classes}

\author{Yan Chen\footnote{Columbia University, Email: yc3107@columbia.edu} ~ and  Jing Dong\footnote{Columbia University, Email: jing.dong@gsb.columbia.edu} \footnote{ Support from NSF grant CMMI-1944209 is gratefully acknowledged by J. Dong} }
\date{}

\maketitle

\begin{abstract}
Utilizing customers' service-time information, we study an easy-to-implement scheduling policy with two
priority classes. By carefully designing the classes, the two-class priority rule achieves near-optimal performance. In particular, 
for a single-server queue, as the traffic intensity approaches 1, the policy achieves a scaling for the queue length processes that is similar to the 
shortest remaining processing time first policy. 
Our analysis quantifies how the tail of the service time distribution affects the benefits one can gain from service-time-based scheduling policies. 
When the service times are misspecified, 
we further quantify how imperfect observation of the service time affects the performance of the two-class priority rule through both theoretical and numerical analysis. Our results demonstrate the robustness of the two-class priority rule. Specifically, even with imperfect service-time information, the two-class priority rules can still achieve substantial performance improvement over the first come first served.  
\end{abstract}




%


\section{Introduction} \label{sec:Intro}
The development of statistical learning techniques and the growing availability of data facilitates us to gain more customer-side information.
From the operations perspective, an important question to ask is how to utilize the customer-side information to achieve better system performance.
In this paper, we look into one particular customer-side information -- service times -- and study how to use service-time information to do smarter scheduling for service systems.

When having perfect service-time information, it is well-known that the shortest remaining processing time first (SRPT) policy achieves superior performance.
In particular, in a single-server queue, it has been established that SRPT is optimal with respect to minimizing the steady-state average sojourn time \cite{SM66} and the average number of customers in the system at each point in time \cite{S68}. 
This policy has been successfully implemented in scheduling jobs for machines (e.g., computers). In those settings, getting job-size information and administering non-first-come-first-served (FCFS) sequencing is relatively easy. Extending this policy to service systems faces two main obstacles: 
(1) managing a queue where the sequence of customers may change upon each arrival is hard, and (2) we may not have perfect job-size information.
Motivated by the main insights from SRPT, in this paper, we propose a two-class priority rule that achieves comparable performance
to SRPT. The two priority classes are defined by a single carefully-chosen threshold where customers whose service times
are below the threshold get higher priority; that is, it prioritizes smaller-size jobs. Within each priority class, customers are served as FCFS. Under this policy, we only need to maintain two FCFS queues,
which greatly facilitate implementations in practice.
As an added benefit of our scheduling policy and analysis framework, we can also provide some theoretical quantification of the effect of imperfect service times. We demonstrate that the policy is robust to service-time mis-specification.

\paragraph{Comparable performance to SRPT.} We start by providing some intuition behind our development.
In queueing systems, smart scheduling is especially important when the system is heavily loaded; that is, when the traffic intensity (server utilization rate), 
$\rho$, is close to $100\%$ (see Section \ref{sec:policy} for a precise definition of $\rho$).
In a single-server queue, as $\rho$ approaches $1$, under FCFS, the queue
scales as $1/(1-\rho)$, whereas under SRPT, the queue can scale much slower than $1/(1-\rho)$ when the service-time distribution has infinite support \cite{LWZ11,P15}.
In other words, SRPT can provide order-of-magnitude performance improvement over FCFS. To see how such improvement is achieved, we notice that if we divide customers into different priority classes according to their remaining service times, i.e., the longer the service time, the lower the priority, then the SRPT policy can be viewed
as the limit of a sequence of multi-class priority rules, where, in the limit, we send the number of classes to infinity \cite{SM66}.
On the other hand, for multi-class queues, under certain conditions on the classes, the properly scaled queue length process only contains jobs in the lowest priority class \cite{R84S}. 
This result is formally known as state-space collapse in heavy-traffic asymptotic analysis. 
Under any work-conserving scheduling policy (the server is not idling when there are still customers waiting in the queue), the total workload scales as $1/(1-\rho)$.
If all the ``remaining" workload is from customers in the lowest priority class, then the queue size scales as $\mu_m/(1-\rho)$, where $\mu_m$ is the service rate (reciprocal of the average service time) of 
customers in the lowest priority class.
The above discussion indicates that if the average service time of the lowest priority class goes to infinity as $\rho$ approaches $1$, then the queue size can scale more slowly than $1/(1-\rho)$.

We make two important observations from the above discussion. First, we only need to divide customers into two classes and focus on how to define the low priority class.
Second, the definition of the lower priority class depends on the service time distribution and the traffic intensity of the system. In particular, as the traffic intensity $\rho$ increases, we may need to shrink the size of the low priority class so that the average service time of jobs in the low priority class increases accordingly. These two observations are the rationale behind the development of our two-class priority rule. The main technical contribution of the paper is to provide a clear rule to define the low-priority class and to rigorously establish that the two-class priority rule with properly defined classes 
achieves comparable scaling for the queue size as SRPT under heavy traffic.


We highlight here the role of service time distribution in how much we can benefit from smart scheduling. 
For our two-class priority rule, the service time distribution determines the appropriate threshold to divide the two classes.
This in turn determines the average service time of the low priority class, and thus the scaling of the queue as $\rho\rightarrow 1$. 
In particular, we show that the value of the threshold depends on the decay rate of the tail cumulative distribution function (tail cdf) of the service time distribution.
The more slowly the tail cdf decays (the heavier the tail of the service time distribution), the more slowly the queue grows as $\rho\rightarrow 1$.
This result indicates that systems in which the service time distributions have heavier tails tend to benefit more from the two-class priority rule.
A similar phenomenon has been observed for the SRPT policy as well \cite{LWZ11}. Our analysis provides a clear explanation of how the service time distribution affects the performance gain from service-time-based priority rules.

\paragraph{Robustness to imperfect service time information.} In real applications, the service-time requirements may not be exactly known to us, and are subject to various estimation errors. There are very limited works studying job size based scheduling with estimation errors, which has been listed as an open problem in \cite{down2019open}. There are several key challenges to study such problems: 1) Estimation errors arise in many different forms, depending on the modeling choice, the estimation error itself may affect the service time information. For example, if the service time is estimated using a linear regression model based on observable characteristics of the customer, then the distribution of the service time is a convolution of the variability in the observables and the estimation error. Thus, different estimation errors may require a very different set of analyses. 2) Process level analysis of system dynamics is challenging as we often have to keep track of either the age process or the remaining processing time process of all the jobs in the system. 3) Optimizing the scheduling policy would require updating the service-time estimation as the job spends more time in service. Our scheduling rule and analysis framework provide an analytically tractable way to investigate how imperfect service-time information affects system performance. In particular, for our two-class priority rule, we only need to study how the estimation error is affecting our ability to classify jobs into the two priority classes. Our asymptotic framework indicates that to achieve order-of-magnitude performance improvement over FCFS, we only require the classified high-priority jobs and low-priority jobs to satisfy certain scaling properties, which can be easy to verify. 

We demonstrate the analytical tractability of our framework by studying models for three simple estimation error models.
Our results generate several insights: 1) when the estimation errors have bounded support, we can achieve the same queue-length scaling as in the perfect information case; 2) when building classification models to assign customers to the two priority classes, it is more important to avoid wrongly classifying short jobs as low-priority jobs than wrongly classifying long jobs as high priority jobs; 3) when dealing with independent and identically distributed (iid) measurement errors, if the tail of the error distribution is lighter than the tail of the service time distribution, we can still achieve $o(1/(1-\rho))$ scaling under the two-class priority rule with a properly chosen threshold. However, when the error distribution has a heavier tail, we may not be able to gain substantial performance improvement from the two-class priority rule. 

We also complement our theoretical analysis with numerical experiments, in which we study the performance of the two-class priority rule with predicted service times. We consider linear regression models, nonlinear regression models, and generalized linear models. In all cases, the two-class priority rule with a properly chosen threshold achieves significant performance improvement over FCFS. Most noticeably, even when the prediction model is only able to explain $3\%$ of the variability in the actual service times, the two-class priority rule utilizing the predicted service times is still able to reduce the steady-state average queue length by more than $50\%$ over FCFS in our simulated example. This is because even with highly inaccurate prediction models, we still have a very low probability of wrongly classifying short jobs as low-priority jobs.

We next provide a brief review of the literature. The goal is to put our work in the right context.

\subsection{Related Literature} \label{sec:lit}

The literature on scheduling using service-time (job-size) information has a long history.
\cite{SM66} are among the first to study the SRPT policy and various other scheduling policies for $M/GI/1$ queues.
Because these scheduling policies require perfect job-size information, we see most applications in scheduling
of jobs for machines \cite{HBSA03}. 
In addition to SRPT, popular scheduling policies include processor sharing (PS), shortest job first. Recently, \cite{SHS18} develop a unified framework to analyze these scheduling policies. 
Relatively few works look at the imperfect information setting \cite{DA14,down2019open}. 
Most noticeably, the foreground-background policy does not require any prior job-size information 
and performs well when the service time distribution has a decreasing failure rate \cite{NW07}.
More generally, policies based on the Gittins index is known to minimize mean response time under various information structures \cite{aalto2009gittins,Scully:2020}. However, the priority structure can be computationally intractable.  
\cite{EIK19} study the benefit of smart scheduling with only predicted service-time information using data from a call center.
Our work contributes to this line of work by proposing a simple two-class priority rule that achieves superior performance and is robust to service-time misspecification.

Our result builds on heavy-traffic asymptotic analysis. Process level asymptotic analysis of job-size-based scheduling policies is challenging, as it often requires keeping track of the remaining processing time of each job in the system. \cite{G04} establishes the diffusion limit for PS queues. \cite{GKP11,P15, banerjee2020heavy} develop the diffusion limit for processes related to SRPT queues. 
Similar to these works, we establish the diffusion limit of the queue length processes under the two-class priority rule. Our scaling is similar to that established in \cite{P15,banerjee2020heavy} for SRPT. 
The heavy-traffic limit of the steady-state average response time under SRPT is studied in \cite{LWZ11}. 
A key observation from this line of work is that the performance of the SRPT policy depends heavily on the tail property of the service time distribution. 
Our work provides a clear explanation for this observation. 

Our work is also related to the line of research on scheduling/prioritizing policies for single-server queues with multiple classes of customers.
\cite{CS61} are among the earliest to analyze the optimality of an important index policy -- the $c\mu$-rule, where $c$ is the per unit time holding cost per customer and $\mu$ is the service rate.
Under the $c\mu$-rule, we give priority to the class with a larger value of  $c\mu$.
Later, \cite{V95} and \cite{MS04} extend the $c\mu$-rule to more general settings using heavy-traffic asymptotic analysis.
In our model, we can think of everyone having the same linear holding cost. Thus, the $c\mu$-rule reduces to one where we prioritize the class with a shorter average service time. 
The main insight from this line of work is that under the diffusion scaling, the queue only contains jobs in the lowest priority class. 
We derive a similar state-space-collapse result. The main difference between our work and the existing literature is that our class division is not fixed. The lower priority class is shrinking as the traffic intensity increases.
The diffusion scaling we apply is also smaller than the one in the $c\mu$-rule literature.
Recent works also study the $c\mu$-rule with imperfect information \cite{AZ09,SAZ18}. Our work also contributes to this line of literature by providing an analytical framework to study the effect of estimation errors. 

Also note that the optimality of scheduling policies with simple structures, such as SRPT or the $c\mu$-rule, has only been established for single-server queues. 
In the asymptotic sense, we can extend the optimality results to multiple-server queues under the conventional heavy-traffic regime where the number of servers is held fixed as the traffic intensity approaches 1. For example, \cite{GZHB18} recently establish the optimality of SRPT for multi-server queues in the conventional heavy-traffic regime. In many-server heavy-traffic regime where the number of servers is sent to infinity with the arrival rate, the asymptotic optimality of $c\mu$-type of scheduling policies may no longer hold (see, e.g., \cite{HZ04,KRW18}).

\subsection{Notations}
Throughout the paper, we refer to waiting time as the total amount of time that a customer spends in the system. This time includes both the time spent in the queue and the time spent being served,
and is also referred to as the sojourn time or the flow time in the literature.
We refer to the number of customers in the queue or the queue size as the number of customers in the system.
This queue length includes both customers waiting to be served and customers in service (being served).
We use the words customer and job interchangeably. 

We define $\eta(t):=t$ and $\zeta(t):=0$ for all $t\geq 0$.
We denote $D[0,\infty)$ as the space of functions from $[0,\infty)$ to $\mathbb{R}$
that are right continuous with left limits, and is endowed with Skorohod $J_1$ topology.
Letting $f(n)$ and $g(n)$ be two nonnegative functions, we write $f(n)\sim g(n)$ if there 
exist $0<c<C<\infty$ and $n_0\in(0,\infty)$, such that for any $n>n_0$, $cg(n)<f(n)<Cg(n)$. We write $f(n)=O(g(n))$ if there exists positive constants $C<\infty$ and $n_0$, such that for any $n \ge n_0$, $f(n) \le C g(n)$.
We write $f(n)=o(g(n))$ if for any $\epsilon>0$, there exists $n_0(\epsilon) \in (0,\infty)$, such that for any $n>n_0(\epsilon)$, $f(n) \le \epsilon g(n)$.

\section{Two-class priority rule} \label{sec:policy}
Consider a single-server queue with a general renewal arrival process and general independent and identically distributed (iid) service times, i.e., a $GI/GI/1$ queue.
We denote $\lambda$ as its arrival rate and $\mu$ as its service rate. 
Then, the traffic intensity $\rho$ is defined as
\[\rho:=\lambda/\mu.\]
We also denote $F$ as the cumulative distribution function (cdf) of the service time, and $\bar F$ as its tail cdf, i.e, $\bar F =1-F$. 
The two priority classes are defined by the threshold
\begin{equation}\label{eq:K}
K(\rho,F)=\bar F^{-1}\left((1-\rho)^{(1-\delta)}\right) \mbox{ for some fixed $0<\delta<1$}.
\end{equation}
Customers whose service times are smaller than or equal to $K(\rho,F)$ are in the high priority class, Class 1,
whereas customers whose service times are larger than $K(\rho, F)$ are in the low priority class, Class 2. $\delta$ is a parameter that balance the load between the high priority class and the low priority class.
We suggest setting $\delta$ small, for example, between $0.01$ and $0.1$, as indicated by our analysis in Section \ref{sec:main}.

The threshold, $K(\rho, F)$, depends on both the service time distribution and the traffic intensity of the system. It is designed
such that $K(\rho, F)$ increases as $\rho$ increases. Note that if the service time distribution has infinite support, 
$K(\rho,F)\rightarrow\infty$ as $\rho\rightarrow 1$. 
However, $K(\rho, F)$ can not grow too fast as $\rho$ approaches $1$. In particular, 
it needs to be properly chosen to control the workload of Class 1 customers (this notion will be made precise in Section \ref{sec:main}).

In Table \ref{tab:per_comp}, we compare the performance of the two-class priority rule with other benchmark policies in $M/M/1$ queues with different traffic intensities.
The benchmark policies include FCFS, shortest-job-first (SJF), where we prioritize the job with a smaller size but do not allow preemption, 
and SRPT, where we prioritize the job with a shorter remaining processing time in a preemptive way.
We consider two implementations of the two-class priority rule, one allows preemption (two-class P) and the other does not (two-class NP).
We observe that compared with FCFS, the two-class priority rule is able to reduce the average queue length by about a half or more across different traffic intensities.
We also note that SRPT performs the best, as has been established in the literature. However, the gaps between the two-class priority rules and SRPT are small.
This finding suggests that even though we can add more priority classes to improve the system performance, the marginal gains from doing so would be small.

\begin{table}[htp]
\centering
\caption{Steady-state average queue length for $M/M/1$ queues under different scheduling policies ($\mu=1$, $\lambda=\rho$, $K=\bar F((1-\rho)^{(1-0.05)})=-(1-0.05)\log(1-\rho)$).} \label{tab:per_comp}
\begin{tabular}{c|cccc}
\hline
$\rho$       & 0.8   & 0.85 & 0.9  & 0.95        \\\hline
FCFS         & 4.00  & 5.67 & 9.00 & 19.00         \\
Two-Class NP & 2.68  & 3.43 & 4.76 & 8.19         \\
Two-Class P  & 2.29  & 2.94 & 4.13 & 7.30   \\
SJF          & 2.31  & 2.86 & 3.78 & 5.95          \\
SRPT         & 1.88  & 2.36 & 3.20 & 5.26   \\ \hline
\end{tabular}
\end{table}

In actual service systems, we may not have perfect service-time information. 
In the next example, we consider the case in which the service times take the form
\[v=\exp(\beta^TX + \epsilon)\]
where $X$ is a vector of observable customer characteristics, $\beta$ is the vector of coefficients for the features in $X$, 
$\epsilon\sim N(0,\sigma_e^2)$ and is independent of $X$.
With this service time model, $\E[\log(v)|X]=\beta^TX$. Assume that we have access to $\beta^TX$. Thus, we can use $\exp(\beta^TX)$
to predict $v$. Let $\hat F$ denote the cdf of the predicted service time, i.e., $\beta^TX$. Customers whose predicted service times are smaller than or equal to $K(\rho, \hat F)$ are assigned to Class 1,
and others are assgined to Class 2.

In Table \ref{tab:per_comp2}, we compare the performance of our two-class priority rule (NP) with FCFS and SJF using the predicted service times.
Set $\beta^TX=0.1+0.1X_1+0.4X_2+0.4X_3$, where $X_1\sim N(1,1)$, $X_2\sim$Exp(1), and $X_3\sim$Uniform$[0,1]$.
We also set $\sigma_e=0.5$. In this case, the prediction model is able to explain 61\% of the variability in the service times.
We keep the service time distribution fixed and vary the arrival rate to achieve different traffic intensities. 
We observe that even with prediction errors, the two-class priority rules can still achieve significant queue length reduction from FCFS. 
In particular, when $\rho=0.95$, the two-class priority rules leads to a more than 50\% queue length reduction over FCFS.
SJF performs slightly better than our proposed policy, suggesting there is limited value in adding more priority classes.
These observations demonstrate the robustness of the job-size based scheduling rules to estimation errors. We provide detailed analysis about how the estimation errors affect the performance of the two-class priority rule in Sections \ref{sec:error} and \ref{sec:num2}.

\begin{table}[htp]
\centering
\caption{Steady-state average queue length for $M/GI/1$ queues under different scheduling policies using predicted service times.} \label{tab:per_comp2}
\begin{tabular}{c|cccc}
\hline
$\rho$       & 0.8   & 0.85 & 0.9  & 0.95        \\
$K(\rho, \hat F)$ & 3.6 & 4 & 4.7 & 6.2 \\\hline
FCFS & 4.54 & 6.60 & 10.04 & 21.21         \\
Two-Class NP  & 3.65  & 4.42 & 6.82 & 10.89   \\
SJF & 2.99  & 3.87 & 5.02 & 8.46   \\\hline
\end{tabular}
\end{table}

In what follows, we start by analyzing the two-class priority rule with perfect service-time information.

\section{Asymptotic performance analysis} \label{sec:main}
To quantify the superior performance of the two-class priority rule and to derive insights on how to choose the appropriate class-division threshold $K(\rho, F)$,
we take a heavy-traffic asymptotic approach. 
We consider a sequence of systems where the traffic intensity converges to $1$ in an appropriate manner,
and we study how the queue length process scales along the sequence.
More specifically, consider a sequence of $GI/GI/1$ queues indexed by $n$, starting from empty at time $0$.
The service time distribution is fixed with mean $1/\mu=1$; that is, we set a unit of time as the average service time. 
The arrival rate for the $n$-th system is $\lambda^n=1-\beta/\sqrt{n}$ for some $\beta>0$.
Under these parameter specifications, the traffic intensity for the sequence of systems approaches 1 at rate $1/\sqrt{n}$, i.e.,  
\[\rho^n:=\lambda^n/\mu=1-\beta/\sqrt{n}.\]

Let $A^n=\{A^n(t): t\geq 0\}$ denote the arrival process of the $n$-th system. 
Let $\tau^n(k)$ denote the interarrival time between the $k$-th arriving customer and the $(k-1)$-th customer. 
We assume $\tau^n(k)$'s are iid and have the same distribution as $\tau^{\infty}/\lambda_n$, where $\tau^{\infty}$ is a random variable with $\E[\tau^{\infty}]=1$. 
We denote $G$ as the cdf of $\tau^{\infty}$. 
We also write $v^n(k)$
as the service time of the $k$-th arriving customer in system $n$.
Because the service time distribution does not change for different scales of systems, we use $v$ to denote a generic 
random variable following the service time distribution $F$.

We impose the following assumption on $v$. 
\begin{assumption}\label{ass:service}
The service time has a continuous distribution with probability density function $f$. There exists $C>0$, such that $f(x)>0$ for any $x>C$. 
There exists $\delta>0$, such that
$\E[v^{2+\delta}]<\infty$.
\end{assumption}
The first part of Assumption \ref{ass:service} essentially requires that the service time distribution has infinite support (unbounded).
Even under SRPT, when the service time distribution has bounded support, the queue 
scales as $1/(1-\rho)$ \cite{LWZ11}. In this case, no order-of-magnitude improvement arises from ``smart" scheduling.
The moment condition in Assumption \ref{ass:service} is standard to establish the diffusion limit for the queue length processes (see, e.g., \cite{WW02}).
Under Assumption \ref{ass:service}, we denote $\sigma_s^2:=\mbox{Var}(v)$.

We also impose the following light-tail assumption on $\tau^{\infty}$.
\begin{assumption}\label{ass:arrive}
There exits $\epsilon>0$ such that $\E[\exp(\theta \tau^{\infty})]<\infty$ for $\theta\in(-\epsilon, \epsilon)$, i.e.,
$\tau^{\infty}$ has a finite moment generating function in a neighborhood of the origin. 
\end{assumption}
Under Assumption \ref{ass:arrive}, $\tau^{\infty}$ has finite moments of all order. We denote $\sigma_a^2:=\mbox{Var}(\tau^{\infty})$.
We also write $\psi_a(\theta): =\log\E[\exp(\theta \tau^{\infty})]$ and $\bar\theta: =\inf\{\theta:\psi_a(\theta)=\infty\}$.
Lastly, we impose the following technical assumption on $\tau^{\infty}$.
\begin{assumption}\label{ass:arrive2}
$\psi_a$ is differentiable everywhere in $(-\infty, \bar\theta)$ and $\lim_{\theta\uparrow\bar\theta}\psi_a(\theta)=\psi_a(\bar\theta)$.
If $\psi_a(\bar \theta)<\infty$, $\lim_{\theta\uparrow\bar\theta}\psi_a^{\prime}(\theta)=\infty$.
\end{assumption}

We note that Assumptions \ref{ass:arrive} and \ref{ass:arrive2} are satisfied by most commonly used interarrival time distributions. 
Examples include phase-type distributions (e.g., Erlang, hyperexponential) and Weibull distributions.

\subsection{The critical threshold value}
For the $n$-th system, we define $K^n$ as the threshold value for the two classes; that is, 
customers with service time less than or equal to $K^n$ are classified into Class 1;
others are classified into Class 2.
For simplicity of exposition, we assume class 1 customers have preemptive priority
over class 2 customers.

We start by introducing a few more notations. 
Denote $\lambda_1^n=\lambda^nF(K^n)$ and $\lambda_2^n=\lambda^n\bar F(K^n)$
as the respective arrival rates of Class 1 and Class 2 customers in the $n$-th system.
We also write  $v_i^n(k)$ as the service time of the $k$-th Class $i$ arrival in system $n$. Then,
\[\mu_1^n=\left(\E[v_1^n(k)]\right)^{-1}=\left(\frac{1}{F(K^n)}\int_{0}^{K^n}xf(x)dx\right)^{-1}\]
and
\[\mu_2^n=\left(\E[v_2^n(k)]\right)^{-1}=\left(\frac{1}{\bar F(K^n)}\int_{K^n}^{\infty} xf(x)dx\right)^{-1}\]
denote the service rates of the two classes respectively.
Let $\rho_i^n=\lambda_i^n/\mu_i^n$ denote the traffic intensity for Class $i$, $i=1,2$.

We next introduce the general idea underlying our development.
From the existing results for priority queues \cite{R84S}, we make the following two important observations:

\begin{enumerate}
\item The Class 1 queue scales as $(1-\rho_1^n)^{-1}$. Note that if Class 1 customers have
preemptive priority over Class 2 customers, the Class 1 queue evolves as if there are only Class 1 customers in the system.
\item The total workload process scales as $(1-\rho^n)^{-1}=O(\sqrt{n})$. 
Indeed, this scaling for the workload process holds for any work-conserving scheduling policies.
\end{enumerate}

Note that if the ``limiting" queue only contains Class 2 customers, the queue length process scales as $\mu_2^n\sqrt{n}$.
To achieve a smaller scaling than $\sqrt{n}$ for the queue length process, we need $\mu_2^n\rightarrow 0$ as $n\rightarrow\infty$. 
However, we need to send $\mu_2^n$ to zero with delicacy.
In particular, we need to make sure that under the scaling $\mu_2^n\sqrt{n}$, we see no Class 1 customers in the queue asymptotically, i.e.,
$\mu_2^n\sqrt{n}(1-\rho_1^n)\rightarrow\infty$ as $n\rightarrow\infty$.

Above all, the key is to choose the threshold $K^n$ in an appropriate way. 
Our choice of $K^n$ satisfies the following assumption.
\begin{assumption} \label{ass:threshold}
$K^n\rightarrow\infty$, and there exists $\delta\in (0,1/2)$, such that
$n^{1/2-\delta}\bar F(K^n)\rightarrow\infty$ 
as $n\rightarrow\infty$.
\end{assumption}
We next provide some comments about Assumption \ref{ass:threshold}.
Let $\gamma^n:=K^n\mu_2^n$.
Under Assumption \ref{ass:service}, $2/3\leq \gamma^n \leq 1$. 
Then, Assumption \ref{ass:threshold} implies that $\mu_2^n\rightarrow 0$
as $n\rightarrow\infty$. This helps ensure that the class 2 queue scales slower than $\sqrt{n}$.
We also note
\[\rho_1^n=\lambda_1^n/\mu_1^n=\left(1-\frac{\beta}{\sqrt{n}}\right)\int_{0}^{K^n}xf(x)dx
=1-\frac{\bar F(K^n)}{\mu_2^n}+ O(1/\sqrt{n}).\]
Thus, under Assumption \ref{ass:threshold}, there exists $\delta\in(0,1/2)$, such that
\begin{equation} \label{eq:rho1}
\frac{n^{1/2-\delta}}{K^n}(1-\rho_1^n)\rightarrow\infty
\mbox{ as $n\rightarrow\infty$,}
\end{equation}
or equivalently,
\[n^{1/2-\delta}\mu_2^n (1-\rho_1^n)\rightarrow\infty
\mbox{ as $n\rightarrow\infty$.}\]
This helps ensure that the class 1 queue scales even slower than the class 2 queue.

Lastly, recall that in \eqref{eq:K}, we suggest setting  $K^n=\bar F^{-1}((1-\rho_n)^{(1-\delta)})$ for some $0<\delta<1$. 
This choice of $K^n$ satisfies Assumption \ref{ass:threshold} when the service time distribution satisfies Assumption \ref{ass:service}.

\subsection{Performance of the two-class priority rule}
In this section, we rigorously quantify the performance of the two-class priority rule with $K^n$ satisfying Assumption \ref{ass:threshold}.
In particular, we analyze how the queue length process scales as $n\rightarrow\infty$.
Denote $Q_i^n=\{Q_i^n(t): t\geq 0\}$ as the queue length process for Class $i$, $i=1,2$, in system $n$,
and $Q^n=Q_1^n+Q_2^n$ as the total queue length. 
The following theorem is the main result of this paper. 
It shows that under the two-class priority rule,
the queue length process 
scales as $\sqrt{n}/K^n$. As $K^n\rightarrow\infty$ as $n\rightarrow\infty$ under assumption \ref{ass:threshold},
the queue length process scales slower than $\sqrt{n}$.

Let $\gamma:=\lim_{n\rightarrow\infty}K^n\mu_2^n$. Define
$\hat Q=\{\hat Q(t): t\geq 0\}$ as a reflected Brownian motion with drift coefficient $-\gamma\beta$
and diffusion coefficient $\gamma\sqrt{\sigma_a^2+\sigma_s^2}$. 

\begin{theorem} \label{th:main}
For the preemptive two-class priority rule, under Assumptions \ref{ass:service} -- \ref{ass:threshold},
\[\frac{K^n}{\sqrt{n}} Q_1^n(nt) \Rightarrow \zeta(t) 
\mbox{ and }
\frac{K^n}{\sqrt{n}} Q_2^n(nt) \Rightarrow \hat Q(t) \mbox{ in $D[0,\infty)$ as $n\rightarrow\infty$},\]
which implies 
\[\frac{K^n}{\sqrt{n}} Q^n(nt) \Rightarrow \hat Q(t)\mbox{ in $D[0,\infty)$ as $n\rightarrow\infty$}.\]
\end{theorem}
Alternatively, we have
\[\frac{1}{\sqrt{n}\mu_2^n} Q^n(nt) \Rightarrow \mbox{RBM}(-\beta, \sigma_a^2+\sigma_s^2) \mbox{ in $D[0,\infty)$ as $n\rightarrow\infty$},\]
where $\mbox{RBM}(-\beta, \sigma_a^2+\sigma_s^2)$ denotes a reflected Brownian motion with drift coefficient $-\beta$
and diffusion coefficient $\sqrt{\sigma_a^2+\sigma_s^2}$. 

Denote $\hat Q(\infty)$ as the stationary distribution of $\hat Q$. Then, $\hat Q(\infty)$ follows an exponential distribution with rate
$2\beta/(\gamma(\sigma_a^2+\sigma_s^2))$ \cite{HW87}.

We also establish the interchange-of-limits and uniform integrability result for the queue length processes of $M/GI/1$ queues. Let $Q^n(\infty)$ denote the stationary distribution of $Q^n(t)$.

\begin{theorem} \label{th:limit_inter}
For the preemptive two-class priority rule, under Assumptions \ref{ass:service} and \ref{ass:threshold}, if the arrival process is Poisson, then for any $x\geq 0$,
\[\lim_{n\rightarrow\infty}\PP\left(\frac{K^n}{\sqrt{n}} Q^n(\infty)\leq x\right)
=1-\exp\left(-\frac{2\beta}{\gamma(1+\sigma_s^2)}x\right)\]
and
\[\lim_{n\rightarrow\infty} \E\left[\frac{K^n}{\sqrt{n}}Q^n(\infty)\right]=\frac{\gamma(1+\sigma_s^2)}{2\beta}.\]
\end{theorem}
Equivalently, we have,
\[\lim_{n\rightarrow\infty} \frac{1}{\sqrt{n}\mu_2^n} \E[Q^n(\infty)] = \frac{1+\sigma_s^2}{2\beta}.\]

  
The proof of Theorem \ref{th:main} is delayed to Section \ref{sec:main_proof},
where we also present several other results about the asymptotic behavior
of the workload processes and the virtual waiting time processes for the two priority classes. 
The proof of Theorem \ref{th:limit_inter} can be found in Appendix \ref{app:thm2}.

\subsection{The effect of the service time distribution} \label{sec:service_tail}
In this section, we study the effect of service time distribution on the performance of the two-class priority rule. We also compare the queue-length scaling under the two-class priority to that under SRPT.
A very appealing aspect of Theorems \ref{th:main} \& \ref{th:limit_inter} is that they explicitly characterize the impact of service time distribution on the benefit of smart scheduling through $K^n$. 
In particular, $Q^n$ scales as $\sqrt{n}/K^n$ and $K^n$ is determined by the service time distribution. In general, the more slowly
the tail of the service time distribution decays, the larger $K^n$ tends to be, which  in turn leads to a smaller scaling of $\E[Q^n(\infty)]$. We demonstrate this through two specific classes of service time distributions. These two classes of distributions include a lot of commonly used service time distributions and are widely studied in the SRPT literature (see, for example, \cite{LWZ11, P15}).

To mark the explicit dependence of threshold $K$ on service time distribution $F$ and the traffic intensity $\rho$, we use the notation $K(\rho, F)$ introduced in \eqref{eq:K}.

\paragraph{Pareto service time distributions}
Consider the family of Pareto distributions: $\bar F_{\alpha}(x)=(m_{\alpha}/x)^{\alpha}$ with parameters $\alpha > 1$ and $m_{\alpha}=(\alpha-1)/\alpha$. In this case,
\[
K(\rho,F_{\alpha})=m_{\alpha}(1-\rho)^{-(1-\delta)/\alpha}.
\]
We note that the smaller the value of the parameter $\alpha$,
the more slowly $\bar F_{\alpha}$ decays, and the larger the value of $K(\rho, F)$ tends to be.

\paragraph{Weibull service time distribution}
Consider the family of Weibull distributions: $\bar F_{\alpha}(x)=\exp(-(x/\nu_{\alpha})^{\alpha})$ with parameters $\alpha>1$ and $v_{\alpha}=1/\Gamma(1+1/\alpha)$, where $\Gamma$ is the gamma function. 
In this case
\[
K(\rho,F_{\alpha})=\nu_{\alpha}\left(-(1-\delta)\log(1-\rho)\right)^{1/\alpha}. 
\]
We again note that the smaller the value of the parameter $\alpha$,
the more slowly $\bar F_{\alpha}$ decays, and the larger the value of $K(\rho, F)$ is.

We next compare the two-class priority rule with SRPT.
For Pareto service time distributions, \cite{LWZ11} establish that under SRPT,
\[n^{-\frac{\alpha-2}{2(\alpha-1)}}\E[Q^n(\infty)]\rightarrow \frac{\pi/(1-\alpha)}{2\sin(\pi/(1-\alpha))}\frac{(1+\sigma_s^2)}{2m_{\alpha}\beta} \mbox{ as $n\rightarrow\infty$.}\]
For our two-class priority rule, if we set $K^n=m_{\alpha} n^{\frac{1}{\alpha}(1/2-\delta)}$ for some $\delta>0$, we have
\[n^{-\frac{\alpha-1-2\delta}{2\alpha}}\E[Q^n(\infty)] \rightarrow \frac{(1+\sigma^2)}{2m_{\alpha}\beta} \mbox{ as $n\rightarrow\infty$.}\]
In this case, even though $\delta$ can be set arbitrarily small, $\frac{\alpha-1}{2\alpha}$ is larger than $\frac{\alpha-2}{2(\alpha-1)}$, i.e., 
the scaling of the average queue length is larger under the two-class priority rule than that under SRPT. However, the difference in scaling decreases as $\alpha$ increases.

For Weibull service time distributions, \cite{LWZ11} establish that under SRPT
$$\frac{\log(\sqrt{n})^{1/\alpha}}{\sqrt{n}}\E[Q^n(\infty)] \rightarrow \frac{1+\sigma_s^2}{2\nu_{\alpha}\beta} \mbox{ as $n\rightarrow\infty$.}$$
For our two-class priority rule, if we set $K^n=\nu_{\alpha}(\log(n^{1/2-\delta}))^{1/\alpha}$ for some $\delta>0$, we have 
\[\frac{\log(n^{1/2-\delta})^{1/\alpha}}{\sqrt{n}}\E[Q^n(\infty)] \rightarrow \frac{1+\sigma_s^2}{2\nu_{\alpha}\beta} \mbox{ as $n\rightarrow\infty$.}\]
In this case, the two scheduling policies achieve the same scaling.

\section{Proof of Theorem \ref{th:main}} \label{sec:main_proof}
In this section, we provide the proof of Theorems \ref{th:main}. 
The proof is divided into two steps.
We first establish state-space collapse results for the workload process and the queue length process
for Class 1 (Proposition \ref{prop:class1}). We then establish proper limits for the 
virtual waiting time process and the queue length process for Class 2 (Proposition \ref{prop:class2}).
The detailed proofs of the propositions are deferred to Appendix \ref{app:proof_prop}.
Note that these intermediate results are of independent interests as they highlight the difference in performance of the two priority classes.

Let $A_i^n=\{A_i^n(t): t\geq 0\}$ denote the arrival process of Class $i$ customers. Note that by our class-division rule, $A_i^n$'s are still renewal processes, but they may not be independent of each other.
We also define
\begin{eqnarray*}
V^n(t)&=&\sum_{k=1}^{A^n(t)}v^n(k)-t,\\
U^n(t)&=&V^n(t)-\inf_{0\leq s\leq t} V^n(s)\wedge 0,\\
V_1^n(t)&=&\sum_{k=1}^{A_1^n(t)}v_1^n(k)-t,\\
U_1^n(t)&=&V_1^n(t)-\inf_{0\leq s\leq t} V_1^n(s)\wedge 0,\\
U_2^n(t)&=&U^n(t)-U_1^n(t).
\end{eqnarray*}
Note that $U^n$ is the total unfinished workload process, 
and $U_i^n$ is the unfinished workload process of Class $i$ jobs.

Under any work-conserving service policy and Assumptions \ref{ass:service} and \ref{ass:arrive}, 
it has been established that \cite{R84S}
\begin{equation}\label{eq:workload}
\frac{1}{\sqrt{n}}U^n(n t) \Rightarrow \mbox{RBM}(-\beta, \sigma_a^2+\sigma_s^2) \mbox{ in $D[0,\infty)$ as $n\rightarrow\infty$}.
\end{equation}

As Class 1 customers have
preemptive priority over Class 2 customers, $U_1^n(t)$ is also the virtual waiting time process for Class 1.
The following proposition establishes state-space collapse for $K^nU_1^n(nt)/\sqrt{n}$
and $K^nQ_1^n(nt)/\sqrt{n}$
\begin{proposition} \label{prop:class1}
Under Assumptions \ref{ass:service} -- \ref{ass:threshold},
$$\frac{K^n}{\sqrt{n}}U_1^n(nt)\Rightarrow \zeta(t) \mbox{ in $D[0,\infty)$ as $n\rightarrow\infty$},$$
and
$$\frac{K^n}{\sqrt{n}}Q_1^n(nt) \Rightarrow \zeta(t) \mbox{ in $D[0,\infty)$ as $n\rightarrow\infty$}.$$
\end{proposition}

\bigskip

Proposition \ref{prop:class1} implies 
the virtual waiting time and the queue length of Class 1 scale as $o(\sqrt{n}/K^n)$, or equivalently, $o(\sqrt{n}\mu_2^n)$. 

Define 
\[\hat U^n(t)=\frac{U^n(nt)}{\sqrt{n}} \mbox{ and }  \hat U_2^n(t)=\frac{U_2^n(nt)}{\sqrt{n}}.\]
From Proposition \ref{prop:class1}, we have
\[\hat U^n-\hat U_2^n\Rightarrow \zeta \mbox{ in $D[0,\infty)$ as $n\rightarrow\infty$,}\]
and thus,
\[\hat U_2^n \Rightarrow \mbox{RBM}(-\beta, (\sigma_a^2+\sigma_s^2)) \mbox{ in $D[0,\infty)$ as $n\rightarrow\infty$.}\]

For Class 2 customers, as it has lower priority, the unfinished workload, $U^n(t)$, is smaller than its virtual waiting time at $t$. 
To characterize the waiting time of Class 2, we define
$$B^n(t,U):=\inf\{s\geq 0: V_1^n(t+s)-V_1^n(t)+U \leq 0\}.$$
Then, $B^n(t,U^n(t))$ is the virtual waiting time process for Class 2.
We also define 
$$\hat B^n(t):=B^n(nt,U^n(nt))/\sqrt{n}.$$

\begin{proposition} \label{prop:class2}
Under Assumptions \ref{ass:service} -- \ref{ass:threshold},
$$\rho_2^n \hat B^n(t) - \hat U^n(t) \Rightarrow \zeta(t) \mbox{ in $D[0,\infty)$ as $n\rightarrow\infty$},$$
and
$$\frac{K^n}{\sqrt{n}}Q_2^n(nt)\Rightarrow \hat Q(t) \mbox{ in $D[0,\infty)$ as $n\rightarrow\infty$}.$$ 
\end{proposition}

Proposition 2 indicates the virtual waiting time process for Class 2 scales as 
\[\sqrt{n}/\rho_2^n \sim \sqrt{n}/(K^n\bar F(K^n)).\]
Note that under FCFS, the virtual waiting time process for Class 2 scales as
\[\sqrt{n}+1/\mu_2^n \sim \sqrt{n}.\] 
Note that under Assumption \ref{ass:threshold}, $K^n\bar F(K^n)\rightarrow0$ as $n\rightarrow\infty$. 
Thus the virtual waiting time for Class 2 under the two-class priority rule has a larger scaling than under FCFS. 
The difference in scaling suggests the overall improvement in average queue length and average waiting time under the two-class priority rule is at the expense of having Class 2 customers waiting for a longer time. 
We refer to \cite{WH03} for more discussions about fairness under different scheduling policies for the $M/GI/1$ queue.

Combining Proposition \ref{prop:class1} and Proposition \ref{prop:class2}, we have proved Theorem \ref{th:main}.

\section{Numerical experiments for the two-class priority rule}
In this section, we provide some numerical experiments to illustrate the pre-limit performance of the two-class priority rule.
In particular, we study the impact of the service time distribution on system performance. We also look into the class-dependent performance of the two-class priority rule.

\begin{figure}[ht]
\caption{The steady-state average queue length for $M/M/1$ queue under the two-class priority rule (P) with different values of the threshold.} \label{fig:threshold_sen}
\centering{
\includegraphics[scale=0.4]{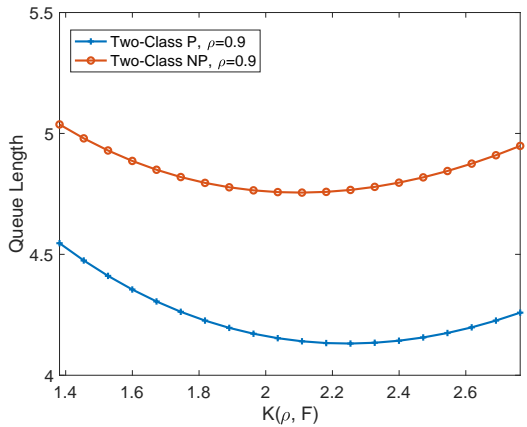} 
}
\end{figure}

We first note that a wide range of $K(\rho,F)$'s satisfies Assumption \ref{ass:threshold}. Thus, we start by conducting some sensitivity analysis for the effect of different values of $K(\rho,F)$ on the performance of the two-class priority rule. 
Figure \ref{fig:threshold_sen} plots the steady-state average queue length for different values of the threshold under the two-class priority rule with and without preemption. Note the expected queue length does vary with different values of $K(\rho,F)$.
However, the magnitude of variation is very small for the range of values of $K(\rho,F)$ plotted.
This finding indicates the two-class priority rule is relatively insensitive to the choice of the threshold within a reasonable range.

In Table \ref{tab:service1}, we compare the steady-state average queue length for $M/GI/1$ queues with the same unit service rate but different service time distributions. 
We consider the class of Pareto distribution with parameter $\alpha$, where we set $m_{\alpha}=(\alpha-1)/\alpha$ and the tail cdf $\bar F(x)=(m_{\alpha}/x)^{\alpha}$. Note that the larger the value of $\alpha$, the faster $\bar F$ decays, and thus, the lighter the tail of the service time distribution is.

\begin{table}[htp]
\centering
\caption{Steady-state average queue length for $M/GI/1$ queues with Pareto($m_{\alpha},\alpha$) service time distributions for different values of $\alpha$ ($\mu=1$, $\lambda=\rho$, $K(\rho,F)=\bar F((1-\rho)^{(1-0.05)})$).} \label{tab:service1}
\begin{tabular}{c|ccc||ccc||ccc}
\hline
& \multicolumn{3}{c||}{FCFS}& \multicolumn{3}{c||}{Two-class (P)} & \multicolumn{3}{c}{SRPT}\\\hline
$\rho$    & 0.8   & 0.9  & 0.99   & 0.8   & 0.9  & 0.99     & 0.8    & 0.9  & 0.99   \\\hline
$\alpha=2.5$& 3.68 &8.19 & 89.20 & 2.53 & 4.43 & 21.60 & 1.88  & 3.05 & 10.54   \\
$\alpha=5$ &  2.51 &5.22 & 53.26  & 2.32 & 4.34 & 30.20 & 2.11  & 3.83 & 23.03  \\
$\alpha=7.5$& 2.44 & 5.05 & 51.18 & 2.38 &  4.57 & 35.91 & 2.20  & 4.16 & 31.64  \\
$\alpha=10$ &  2.42 & 5.00 & 50.61 & 2.40  & 4.70 & 39.16 & 2.25   & 4.35 & 34.81 \\\hline
\end{tabular}
\end{table}

We observe that when comparing the two-class priority rule to FCFS, the two class-priority rule always achieves a smaller queue. The queue length reduction under the two-class priority rule is becoming larger when the traffic intensity increases. For example for Pareto service distribution with $\alpha=2.5$, when $\rho=0.8$, the two-class priority reduces the average queue length by $32\%$; when $\rho=0.9$, the reduction is $46\%$; in the extreme case when $\rho=0.99$, the reduction is as large as $76\%$. These observations are consistent with our asymptotic analysis in Theorems \ref{th:main} \& \ref{th:limit_inter}. Specifically, the queue length process scales slower as $\rho$ increases under the two-class priority rule than that under FCFS, and this scaling difference is more apparent when $\rho$ is closer to 1.

We next take a closer look at the effect of the service time distributions.
When comparing the two-class priority rule with FCFS, we note that we gain more reduction in queue length for heavier tail service time distributions, i.e., when $\alpha$ is smaller. For example, for $\rho=0.9$, when $\alpha=2.5$, we achieve a queue length reduction of $46\%$ using the two-class priority rule; when $\alpha=5$, the reduction is $17\%$; when $\alpha=10$, the reduction is only $6\%$.

We also oberve that under FCFS, as the tail of the service time distribution becomes lighter, i.e., $\alpha$ increases, the average queue length decreases. This finding is expected, because for $M/GI/1$ queue with unit service rate, the steady-state expected queue length under FCFS takes the form
\[
\rho+\frac{(1+\sigma_s^2)\rho^2}{2(1-\rho)},
\]
and $\sigma_s^2$ decreases as $\alpha$ increases. 
Under the two-class priority rule, Theorem \ref{th:limit_inter} suggests that 
\[\E[Q(\infty)]\approx \frac{\gamma(1+\sigma_s^2)}{2(1-\rho)K(\rho, F_{\alpha})} \mbox{ where }K(\rho,F_{\alpha})=m_{\alpha}(1-\rho)^{-(1-\delta)/\alpha}.\]
Note that when $\alpha$ increases, both $\sigma_s^2$ and $K(\rho, F_{\alpha})$ decreases. 
As $\sigma_s^2$ does change with $\rho$, the effect of $K(\rho, F_{\alpha})$ will  
dominate when $\rho$ is close enough to 1. 
However, when $\rho$ is small, whether the smaller $\sigma_s^2$ or the smaller $K^n$ will dominate is unclear. 
From Table \ref{tab:service1}, we observe that when $\rho=0.8$, the average queue length is first decreasing and then increasing as $\alpha$ increases. In this regime, the effect of $\sigma_s^2$ plays a role. However, when $\rho=0.99$, the average queue length increases as $\alpha$ increases. In this regime, the effect of $K(\rho, F_{\alpha})$ dominates.

Lastly, when comparing SRPT with the two-class priority rule, we note that when $\alpha$ is small, 
SRPT can achieve substantially shorter average queue length than the two-class priority rule. However,
the expected queue lengths under the two scheduling policies are getting closer as $\alpha$ increases.
For example, for $\rho=0.9$, when $\alpha=5$, the optimality gap of the two-class priority rule is $13\%$;
when $\alpha=10$, the gap is $8\%$. This observation is consistent with our analysis in Section \ref{sec:service_tail}.

\begin{figure}[ht]
\caption{Steady-state average queue lengths of $M/M/1$ queues with different traffic intensities under the two-class priority rule (P).} \label{fig:queue}
\centering{
\includegraphics[scale=0.4]{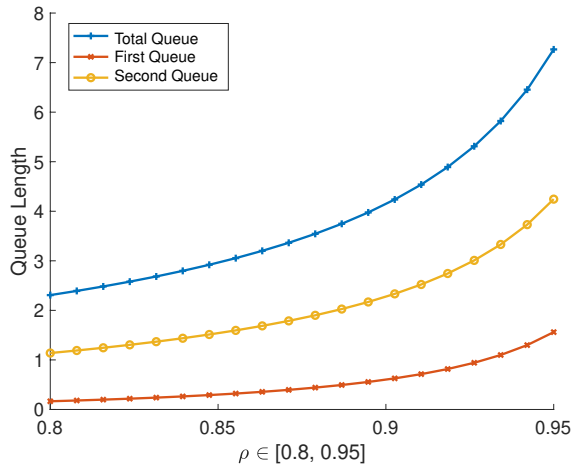} 
}
\end{figure}

We next demonstrate the class-dependent performance under the two-class priority rule.
Figure \ref{fig:queue} compares the steady-state average queue length for Class 1 and Class 2 under the two-class priority rule for an $M/M/1$ queue. We also plot the average total queue length. We observe that because Class 1 is prioritized, it has a much shorter queue than Class 2. Table \ref{tab:wait} compares the steady-state average waiting times for Class 1, Class 2, and all customers together (All) under the two-class priority rule and FCFS. We observe that compared with FCFS, despite the huge improvement in average waiting times for Class 1 customers, Class 2 customers do incur a significant increase in waiting times under the two-class priority rule. These observations are consistent with our asymptotic analysis in Propositions \ref{prop:class1} \& \ref{prop:class2}.

\begin{table}[htp]
\centering
\caption{Steady-state average waiting time for $M/M/1$ queues under the two-class priority rule (P) and FCFS ($\mu=1$, $\lambda=\rho$, $K(\rho)=\bar F((1-\rho)^{(1-0.05)})$).} \label{tab:wait}
\begin{tabular}{c|ccc||ccc}
\hline
&  \multicolumn{3}{c||}{Two-Class (P)} & \multicolumn{3}{c}{FCFS}\\\hline
 & Class 1   & Class 2  & All   & Class 1   & Class 2 & All \\\hline
$\rho=0.8$ & 0.83 & 10.19 & 2.86 & 4.58 & 6.53 & 5.00  \\
$\rho=0.85$ &  1.07 & 15.61 & 3.46  & 6.31 & 8.47 & 6.67   \\
$\rho=0.9$ & 1.52 & 28.98 & 4.59 & 9.72 &  12.19 & 10.00   \\
$\rho=0.95$ &  2.79 & 89.34 & 7.68 & 19.82  & 22.85 & 20.00  \\\hline
\end{tabular}
\end{table}

\section{Imperfect service-time information} \label{sec:error}

As an added benefit to our two-class priority rule, we can provide some theoretical quantification for the effect of imperfect service-time information on system performance. Employing the same heavy-traffic asymptotic mode of analysis as in Section \ref{sec:main},
our goal is to be able to gain ``order-of-magnitude" performance improvement using the two-class priority rule over FCFS. When we do not have perfect job-size information, what we classified as Class 1 jobs may be different from jobs whose actual service times are less than or equal to $K^n$. For $i=1,2$, let $\hat\lambda_i^n$, $\hat\mu_i^n$, and $\hat\rho_i^n$ denote the arrival rate, service rate, and traffic intensity of the classified Class $i$ jobs respectively.
The key insights from our perfect information analysis are that as long as the following two conditions hold, the queue length under the two-class priority rule still scales slower than $(1-\rho_n)^{-1}$.
\begin{itemize}
\item[C1)] The average service time of the classified Class 2 jobs goes to infinity as the traffic intensity approaches 1. In particular, if $\hat \mu_2^n\rightarrow 0$ as $n\rightarrow\infty$, the Class 2 queue still scales slower than $(1-\rho_n)^{-1}$.
\item[C2)] The queue of the classified Class 1 jobs under a proper scaling (determined by the size of the Class 2 queue) diminishes as the traffic intensity approaches 1. In particular, if $n^{1/2-\delta}\hat\mu_2^n(1-\hat\rho_1^n)\rightarrow\infty$ as $n\rightarrow\infty$ for some $0<\delta<1/2$, the Class 1 queue still scales slower than the Class 2 queue.
\end{itemize}
These two conditions provide quite some flexibility to accommodate estimation errors for the service times. 
In our subsequent development, we study how different estimation errors would affect the above two conditions (C1 and C2).

We denote $v^n(k)$ as the actual size of the $k$-th job in the $n$-th system and $\hat v^n(k)$ as its estimated size. We also write $\varepsilon^n(k)=v^n(k)-\hat v^n(k)$ as the estimation error. As the service time does not scale with $n$, we denote $v$ as a generic service time, $\hat v$ as a generic estimated service time, and $\varepsilon$ as a generic estimation error.  We refer to the real Class 1 jobs as jobs whose actual service times are no larger than the threshold $K^n$, i.e, $v^n(k)\leq K^n$, and similarly for the real Class 2 jobs. We refer to the classified Class 1 jobs as jobs who are assigned to Class 1 based on some estimation, and similarly for the classified Class 2 jobs. 

We consider three specific forms of estimation errors. 
The first one has bounded estimation errors. We employ a worst-case analysis to demonstrate the robustness of the two-class priority rule.
The second one is a simple classification model, which allows us to investigate the impact of two different types of classification errors: wrongly classifying a Class 1 job as Class 2 versus wrongly classifying a Class 2 job as Class 1. 
Our third model assumes iid measurement errors, which allows us to investigate the interplay between the tail of the measurement error distribution and the tail of the service time distribution. Even though these are simplified forms of estimation errors, they provide important 
insights into the effect of estimation errors on system performance under the two-class priority rule.
In Section \ref{sec:num2}, we complement the theoretical analysis in this section with numerical experiments for problems with more general predicted service times and in more realistic settings.

\subsection{Bounded estimation errors} \label{sec:bounded}
Consider the case where there exits $M\in(0,\infty)$ such that $\varepsilon^n(k) \in [-M, M]$ with probability one for all $k\in \mathbb{N}$.
We do not impose further assumptions on the estimation errors. For example, $\varepsilon^n(k)$ can have non-zero mean and can depend on $v^n(k)$ and/or $\hat v^n(k)$. For the two-class priority rule, we assign the jobs whose predicted service time is smaller than or equal to $K^n$, i.e., $\hat v^n(k)\leq K^n$, to Class 1; and the others are assigned to Class 2. 

With bounded estimation errors, we can take a worst case analysis approach. In particular, the worst-case scenario that is ``against" our two-class priority rule is the one where we wrongly classify all jobs whose actual service times are in $[K^n-M, K^n]$ as Class 2, and all the jobs whose actual service times are in $(K^n, K^n+M]$ as Class 1. (This notion will be made precise in the proof of Proposition \ref{th:bounded}.) Even in this worst-case scenario, we are still able to show that $1/\hat \mu_2^n \sim K^n$ and $n^{1/2-\delta}\hat \mu_2^n(1-\rho_2^n)\rightarrow\infty$ as $n\rightarrow\infty$ for some $0<\delta<1/2$. Thus, we achieve the same scaling for the queueing length process as in the perfect information setting:  

\begin{proposition}\label{th:bounded}
For the preemptive two-class priority rule with estimated service times, under Assumptions \ref{ass:service} -- \ref{ass:threshold}, suppose there exits a constant $M\in(0,\infty)$ such that $\varepsilon \in [-M, M]$ with probability one. Then, $1/\hat \mu_2^n \sim K^n$ and
\[
\frac{1}{\sqrt{n}\hat \mu_2^n} Q^n(nt) \Rightarrow \mbox{RBM}(-\beta, \sigma_a^2+\sigma_s^2)  \mbox{ in $D[0,\infty)$ as $n\rightarrow\infty$}.
\]
\end{proposition}

The proof of Proposition \ref{th:bounded} and all subsequent results in this section can be found in Appendix \ref{app:error}.


\subsection{A special classification model} \label{sec:class}
In this section, we assume that there exists a classification model that classifies incoming jobs into the two priority classes.
The classification errors depend on the underlying classification model. For example, if there is a prediction model for service times, e.g, the model considered in Section \ref{sec:bounded}, and we classify  jobs whose predicted service times are less than or equal to $K^n$ as Class 1, then the probability of wrongly classifying a Class 2 job as Class 1 depends on the actual service time of the job and the accuracy of the prediction model. 
In this case, jobs whose actual service times are closer to $K^n$ may be more likely to be wrongly classified.

In this section, we consider a simplified setting where for the $n$-th system, each real Class 1 job has an equal probability of being wrongly classified as a Class 2 job. We denote this probability as $p_{12}^n$. Likewise, each real Class 2 job has an equal probability of being wrongly classified as a Class 1 job. We denote this probability as $p_{21}^n$. The goal is to highlight the difference between the two types of classification errors, i.e, $p_{12}^n$ versus $p_{21}^n$.
 
Let $p_{11}^n=1-p_{12}^n$ and $p_{22}^n=1-p_{21}^n$. Then, we have
\[\begin{split}
\hat \lambda_i^n &= p_{1i}^n\lambda^n F(K^n) + p_{2i}^n\lambda^n\bar F(K^n),\\
\hat \mu_i^{n} &=\left(\frac{p_{1i}^n \lambda^n F(K^n)}{\hat \lambda_i^n} \frac{\int_{0}^{K^n}xf(x)dx}{F(K^n)}    
+\frac{p_{2i}^n \lambda^n \bar F(K^n)}{\hat \lambda_i^n} \frac{\int_{K^n}^{\infty} xf(x)dx}{\bar F(K^n)}\right)^{-1}.
\end{split}\]

\begin{proposition} \label{th:class}
For the preemptive two-class priority rule with a classification model where each real Class $i$ job has equal probability of being classified as Class $j$, $i,j=1,2$, assume Assumptions \ref{ass:service} -- \ref{ass:threshold} hold.
\begin{itemize}
\item[i)] If $p_{12}^n=O(\bar F(K^n))$ and $p_{21}^n\leq a$ for some $0\leq a<1$, then $1/\hat \mu_2^n \sim K^n$ and
\[
\frac{1}{\sqrt{n}\hat \mu_2^n} Q^n(nt) \Rightarrow \mbox{RBM}(-\beta, \sigma_a^2+\sigma_s^2)  \mbox{ in $D[0,\infty)$ as $n\rightarrow\infty$}.
\]
\item[ii)] If $p_{12}^n/(K^n\bar F(K^n))\rightarrow\infty$ as $n\rightarrow\infty$ and $p_{21}^n\leq a$ for some $0\leq a<1$, then $\lim_{n\rightarrow\infty}\hat \mu_2^n = 1$.
\end{itemize}
\end{proposition}

Proposition \ref{th:class} indicates that it is more important to avoid wrongly classifying real Class 1 customers as Class 2 than the other way around, i.e., it is more important to avoid wrongly classify short jobs as long ones than the other way around.
In particular, as long as $p_{12}^n$, decays to zero sufficiently fast and the probability of the second type of error, $p_{21}^n$, is bounded away from $1$, we can achieve the same queue-length scaling as in the perfect information case. On the other hand, if $p_{12}^n$ decays to 0 too slowly, the queue may scale the same as that under FCFS (Case ii in Proposition \ref{th:class}). This observation can be counterintuitive at first glance, but note that to gain ``order-of-magnitude" performance improvement under the two-class priority rule, we need $\hat\mu_2^n\rightarrow 0$ as $n\rightarrow \infty$. If we wrongly classify too many short jobs as long ones, we will not be able to achieve this convergence. On the other hand, the classified short jobs are being processed in an ``underloaded" regime, thus their queue is small.

\subsection{Measurement errors} \label{sec:measure}
To gain more insights into the interplay between the service time distribution and the distribution of estimation errors, in this section, we consider a specific form of estimated service times that may arise due to measurement errors. We assume
$\varepsilon^n(k)$'s are iid with finite mean and variance, and $\varepsilon^n(k)$ is independent of $v^n(k)$.
Under the two class priority rule, if the estimated job size is smaller than or equal to $K^n$, we assign the job to Class 1; 
otherwise, the job is assigned to Class 2.

In this case, without loss of generality, we write 
\[\hat v^n(k) = v^n(k) + \varepsilon^n(k)\]
instead of $\hat v^n(k) = v^n(k) - \varepsilon^n(k)$.
We assume $\varepsilon$ is a continuous random variable with probability density function $\phi$. We also write $\Phi$ as the cdf of $\varepsilon$ and $\bar \Phi:=1-\Phi$ as its tail cdf.
Then, 
\begin{eqnarray*}
&&\hat \lambda_1^n=\lambda^n \PP(v+\varepsilon<K^n)=\left(1-\frac{\beta}{\sqrt{n}}\right)\int_0^\infty f(t)\Phi(K^n-t)dt;\\
&&\hat \lambda_2^n=\lambda^n \PP(v+\varepsilon \ge K^n)=\left(1-\frac{\beta}{\sqrt{n}}\right)\int_0^\infty f(t)\bar \Phi(K^n-t)dt;\\
&&1/\hat \mu_1^n=\E[v|v+\varepsilon<K^n]=\frac{\int_0^\infty tf(t)\Phi(K^n-t)dt}{\int_0^\infty f(t)\Phi(K^n-t)dt};\\
&&1/\hat \mu_2^n=\E[v|v+\varepsilon \ge K^n]=\frac{\int_0^\infty tf(t)\bar\Phi(K^n-t)dt}{\int_0^\infty f(t)\bar\Phi(K^n-t)dt}.
\end{eqnarray*}


We first note that
\begin{eqnarray*}
 1-\hat \rho_1^n & = & 1-\left(1-\frac{\beta}{\sqrt{n}}\right)\int_0^\infty tf(t)\Phi(K^n-t)dt \\
& = & \int_0^\infty tf(t)\bar \Phi(K^n-t)dt +O(1/\sqrt{n}).
\end{eqnarray*}
In addition, note that for $\delta>0$ chosen according to Assumption \ref{ass:threshold},
\[\begin{split}
n^{1/2-\delta}\hat\mu_2^n(1-\hat \rho_1^n)&= n^{1/2-\delta}\int_0^\infty f(t)\bar \Phi(K^n-t)dt + o(1)\\
&\geq n^{1/2-\delta} \int_{K^n}^\infty f(t)\bar \Phi(K^n-t)dt+o(1)\\
&\geq \bar\Phi(0) n^{1/2-\delta}\bar F(K^n) \rightarrow \infty \mbox{ as $n\rightarrow\infty$.}
\end{split}\]
This indicates that $n^{1/2-\delta}\mu_2^n(1-\rho_n) \rightarrow\infty$  as $n\rightarrow\infty$.
In particular, even with estimation errors, the Class 1 queue still scales slower than the Class 2 queue.
Thus, our main focus in subsequent analysis is to decide whether $\hat \mu_2^n\rightarrow 0$ as $n\rightarrow\infty$,
i.e., to see whether the Class 2 queue scales slower than $(1-\rho^n)^{-1}$.

The following assumption put restrictions on the tail of the distribution of $\varepsilon$. In particular, it requires that $\varepsilon$ 
has a lighter tail than $v$.
\begin{assumption}\label{ass:error2}
There exists $c>0$ such that for any $0<\theta<c$, $\E[\exp(-\theta\varepsilon)]<\infty$. 
For any $M>0$, $\lim_{t\rightarrow\infty} \bar \Phi(t-M)/\bar F(t)=0$.
\end{assumption}

%
%
%
%

\begin{proposition}\label{th:est} 
For the preemptive two-class priority rule using estimated service times with iid measurement errors, under Assumptions \ref{ass:service} - \ref{ass:error2},
$\hat \mu_2^n \rightarrow 0$ and
\[\frac{1}{\sqrt{n}\hat \mu_2^n} Q^n(nt) \Rightarrow \mbox{RBM}(-\beta, \sigma_a^2+\sigma_s^2)  \mbox{ in $D[0,\infty)$ as $n\rightarrow\infty$}.\]
If we further assume that for any $a>0$, $\lim_{t\rightarrow\infty} \bar \Phi(at)/\bar F(t)=0$, then
$1/\hat \mu_2^n \sim K^n$.
\end{proposition}

Proposition \ref{th:est} indicates that if the estimation error distribution has a lighter tail than the service time distribution, we can still achieve $o(\sqrt{n})$ scaling for the queue length process in heavy traffic. 
From the proof of Proposition \ref{th:est}, we note that the exact order of the scaling may depend on ``how much lighter" the estimation error distribution is in comparison to the service time distribution. 


Given Theorem \ref{th:est}, a natural follow-up question is what would happen if the estimation error  distribution has a heavier tail than the service time distribution.
We next study a special case of this:

\begin{lemma}\label{lem:heavy}
Suppose the distribution of the estimation error is regularly varying. In particular, there exists a function $\xi:(0,\infty) \rightarrow(0,\infty)$ such that for any $a>0$, 
$\lim_{x\rightarrow\infty} \bar \Phi(ax)/\Phi(x)=\xi(a)$.
Under Assumptions \ref{ass:service} - \ref{ass:threshold}, 
if
\[\lim_{x\rightarrow\infty} \frac{\bar \Phi(x)}{x\bar F(x)}=c<\infty,\]
then $\lim_{n\rightarrow\infty} \hat \mu_2^n>0$.
\end{lemma}

Lemma \ref{lem:heavy} requires the distribution of $\varepsilon$ to be regularly varying, For example, the error distribution can have a Pareto tail, 
$\bar\Phi(x)=\left(\frac{m}{x+m\alpha/(\alpha-1)}\right)^{\alpha}$. It also requires that the error distribution has a heavier tail than the service time distribution.
For example, if the service time distribution is Pareto with parameter $\alpha_s$, $\alpha_s>3$, and the error distribution has a Pareto tail with parameter $\alpha_e\in(1,\alpha_s-1]$,
the conditions of Lemma \ref{lem:heavy} are satisfied. In this case, the two-class priority rule gives rise to the same queue-length scaling as FCFS.




\section{Numerical experiments for predicted service times} \label{sec:num2}
In this section, we consider several classical prediction models for the service times, and
use simulation to evaluate the performance of the two-class priority rule when we only have access to predicted service times.
Focusing on more realistic service system settings, we study the non-preemptive version of the two-class priority rule (2NP), and choose FCFS and SJF with predicted service times as two benchmark policies. All the reported steady-state average queue length is estimated based on the average queue length over $10^7/\mu$ units of time (long-time average). To facilitate comparison, common random numbers are used when simulating the system under different scheduling policies.


\begin{example}[A linear regression model] \label{eg1}
Consider the setting where we know certain features of the customers and the service time model takes the form
\[v_i=\beta^TX_i + \epsilon_i,\]
where $X_i$ denotes the vector of observable features of the $i$-th customer, $\beta$ is the vector of coefficients, and $\epsilon_i$ is the residual that is assumed to have mean zero and is uncorrelated with $X_i$.
Note that in this case, $\E[v_i|X_i]=\beta^TX_i$.
\end{example}

In Table \ref{tab:eg1}, we consider the linear model:
$\E[v_i|X_i]=0.1+0.1X_{1i} + 0.4X_{2i} + 0.4 X_{3i}$, where $X_{1i}\sim N(1,1)$, $X_{2i} \sim$Exponential(1), and $X_{3i} \sim$Uniform$[0,2]$, 
and $\epsilon_i$'s are iid $N(0,\sigma_e^2)$ and are independent of $X_i$'s.
Note that in this case, $\epsilon_i$ and $v_i$ are dependent. Thus, this model is different from the measurement error model studied in Section \ref{sec:measure}.
When generating service times, we set 
\begin{equation}\label{eq:trunc}
v_i=(\E[v_i|X_i]+\epsilon_i)\vee 0.001
\end{equation}
to ensure that the service times are positive. 
The value of $\sigma_e$ varies form $0$ to $0.5$. 
Due to \eqref{eq:trunc}, 
setting $\sigma_e$ too large would require a non-negligible amount of truncations which makes the service time model quite different from the linear regression model. Thus, $\sigma_e$ is capped at $0.5$. 
When running 2NP and SJF, we use $\beta^TX_i$ as the predicted service time.
The traffic intensity of the system is kept at $95\%$ and the threshold of the two-class priority rule is set to be the 95-th quantile of the predicted service time distribution, which is equal to $1.8$ in this case.

We observe from Table \ref{tab:eg1} that as $\sigma_e$ increases from $0$ to $0.5$, the queue length reduction we gain from 2NP over FCFS only decreases slightly, i.e., from $36\%$ when $\sigma_e=0$ to $35\%$ when $\sigma_e=0.5$. This suggests that 2NP is quite robust to estimation errors.
We also observe that the two-class priority rule achieves comparable but worse performance than SJF. 

\begin{table}[htp]
\centering
\caption{Steady-state average queue length of $M/G/1$ queues with predicted service times $\beta^TX_i$. $\lambda=0.95\mu$. For the two class priority rule, we set the threshold $K=1.8$.
} \label{tab:eg1}
\begin{tabular}{c|ccc}
\hline
$\sigma_e$       & 0   & 0.25  & 0.5        \\\hline
2NP & 7.67 & 7.92 & 8.31        \\
SJF  & 6.39  & 6.50 & 6.96   \\
FCFS & 12.00  & 12.35 & 12.79     \\\hline
\end{tabular}
\end{table}

\begin{example}[A log-linear regression model]
Consider the following service time model
\[v_i=\exp(\beta^TX_i + \epsilon_i),\]
where $X_i$ denotes the vector of observable features for the $i$-th customer, $\beta$ is the vector of coefficients, 
and $\epsilon_i$ has mean zero and is uncorrelated with $X_i$. 
In this case, $\E[\log(v_i)|X_i]=\beta^T X_i$.
\end{example}

In Table \ref{tab:eg2}, we set
$\E[\log(v_i)|X_i]=0.1+0.1X_{1i} + 0.4X_{2i} + 0.4 X_{3i}$, where $X_{1i}\sim N(1,1)$, $X_{2i} \sim$Exponential(1), and $X_{3i} \sim$Uniform$[0,2]$, 
and $\epsilon_i$'s are iid $N(0,\sigma_e^2)$ and are independent of $X_i$'s.
The value of $\sigma_e$ is varied from $0$ to $1.5$. 
When running 2NP and SJF, we use $\exp(\beta^TX_i)$ as the predicted service time. Note that in this case, the predicted service time is biased, i.e., it does not have the same expectation as the actual service time, and the estimation errors conditional on the observable features, i.e, $v_i-\exp(\beta^TX_i)$ conditional on $X_i$, are not identically distributed. The traffic intensity of the system is kept at $95\%$ and the threshold of the two-class priority rule is set to be the 95\% quantile of the predicted service time distribution, which is equal to $6.2$ in this case.

We first observe that across all levels of estimation errors tested in Table \ref{tab:eg2}, 2NP leads to a more than 50\% reduction in the steady-state average queue length when comparing to FCFS. 
In the extreme case when $\sigma_e = 1.5$, the predicted service times explain only $3\%$ of the variability in the actual service times. Even in this case, 2NP achieves a 51\% reduction in average queue length. We also note that under 2NP, the proportion reduction is not monotonically decreasing as $\sigma_e$ increases. This is because in this example, as $\sigma_e$ increases, the coefficient of variation of the service time distribution increases significantly. For example, when $\sigma_e=0$, the coefficient of variation is only $0.92$; when $\sigma_e=1.5$, the coefficient of variation is $4.50$. When the service times are more heterogeneous, we may gain more from prioritizing shorter jobs. Lastly, we observe that SJF achieves an even smaller average queue length than 2NP. For example, when $\sigma_e=0.5$, 2NP leads to a $49\%$ reduction while SJF leads to a 60\% reduction in average queue length over FCFS. When $\sigma_e=1$, 2NP leads to a $55\%$ reduction while SJF leads to a 68\% reduction. 
This suggests that we can gain further performance improvement by introducing more priority classes, even with estimation errors. Note that SJF can be viewed as the limit of a properly defined sequence of multi-class priority rules as the number of priority classes goes to infinity \cite{SM66}.  
 
\begin{table}[htp]
\centering
\caption{Steady-state average queue length of $M/G/1$ queues with predicted service times $\exp(\beta^TX_i)$. $\lambda=0.95 \mu$. For the two class priority rule, we set the threshold $K=6.5$.
} \label{tab:eg2}
\begin{tabular}{c|ccccccc}
\hline
$\sigma_e$    & 0   & 0.25 & 0.5  & 0.75  & 1 & 1.25   & 1.5  \\\hline
2NP & 8.36 & 9.20 & 10.89 & 14.32    &  20.76 & 34.70 & 67.58\\
SJF  & 6.52  & 6.74 & 8.46 & 10.09  & 14.75 & 24.06 & 47.18 \\
FCFS & 16.82  & 18.47 & 21.21 & 30.67 & 46.00 & 74.52 & 137.54  \\\hline
\end{tabular}
\end{table}

We next take a closer look at why even with very noisy predictions, utilizing the predicted service times to do smart scheduling can still achieve significant performance improvement. Table \ref{tab:eg22} reports the $R^2$, which measures the proportion of variability in the actual service times that can be explained by the prediction model, $p_{12}$, which is the probability of wrongly classifying a Class 1 job into Class 2, and $p_{21}$, which is the probability of wrongly classifying a Class 2 job into Class 1. We observe that even though $R^2$ decreases substantially and $p_{21}$ increases substantially as $\sigma_e$ increases, $p_{12}$ only increases by a very small amount and is still very close to zero. In particular, $p_{12}$ increases from $0.011$ when $\sigma_e=0.25$ to $0.033$ when $\sigma_e=1.5$. Based on our analysis in Section \ref{sec:class}, $p_{12}$ plays a more important role in determining the performance of the two-class priority rule. In particular, Proposition \ref{th:class} indicates as long as $p_{12}$ is sufficiently small and $p_{21}<1$, we can achieve similar performance as in the perfect information case. 

\begin{table}[htp]
\centering
\caption{Estimation errors for Example 2} \label{tab:eg22}
\begin{tabular}{c|cccccc}
\hline
$\sigma_e$   & 0.25 & 0.5  & 0.75  & 1 & 1.25   & 1.5  \\\hline
$R^2$       &  87\% & 61\% & 33\% & 20\% & 8\% & 3\% \\
$p_{12}$   & 0.011 & 0.017 & 0.022  & 0.026 & 0.029 & 0.033 \\
$p_{21}$   & 0.324 & 0.630 & 0.778 & 0.837 & 0.868 & 0.886  \\\hline
\end{tabular}
\end{table}


\begin{example}[A Gamma model]
Consider a simplified supermarket checkout setting where we get to observe the number of items a customer is purchasing, which we denoted by $L_i$ for the $i$-th customer.
Given $L_i$, we assume the actual service time $v_i$ has a Gamma distribution with parameters $L_i$ and $\theta$, where $\theta$ is average processing time of each item, which we normalize to be one unit of time. 
\end{example}

In Table \ref{tab:eg3}, we assume $L_i$'s are iid Poisson random variables with rate $10$. When running 2NP and SJF, we use $L_i$ as the predicted service time. The traffic intensity $\rho$ varies from $80\%$ to $95\%$. For a particular value of $\rho$, the threshold $K$ is set to be the $100\rho\%$ quantile of $L_i$. 

We observe from Table \ref{tab:eg3} that as $\rho$ increases, we gain more performance improvement from 2NP over FCFS. Specifically, when $\rho=0.8$, the reduction in steady-state average queue length under 2NP is only $13\%$. When $\rho=0.95$, the reduction is $22\%$. In addition, we observe again that SJF achieves a slighter better performance than 2NP. However, SJF may be harder to implement in practice as the sequence of waiting customers may change upon each new arrival.

\begin{table}[htp]
\centering
\caption{Steady-state average queue length of $M/G/1$ queues with predicted service times $L_i$. $\lambda=0.95/10$. 
} \label{tab:eg3}
\begin{tabular}{c|cccc}
\hline
$\rho$ & 0.8   & 0.85  & 9 & 0.95        \\
$K$ & 12 & 13 & 14 & 15 \\\hline
2NP & 2.36 & 3.20 &  4.79 & 9.15        \\
SJF  & 2.27  & 3.00 & 4.32 & 8.08   \\
FCFS & 2.71  & 3.69 & 5.72 & 11.72     \\\hline
\end{tabular}
\end{table}

\section{Concluding remarks}
In this paper, we propose a simple two-class priority rule that achieves comparable performance to SRPT.
We characterize the process-level diffusion limit for the queue length process in heavy-traffic.
The diffusion scaling is nonstandard and depends on the service time distribution. 
We also demonstrate the robustness of the two-class priority rule to service-time misspecification through both theoretical analyses
and numerical experiments.

Our theoretical analysis of service-time misspecification focuses on simplified models of estimation errors to gain analytical tractability. Our numerical experiments consider more general regression models for service times and demonstrate the good performance of the two-class priority rule. Theoretical analysis of the two-class priority rule with more general predicted service times would be an interesting future research direction.
We also observe in our numerical experiments that SJF with predicted service times, in general, achieves better performance than the two-class priority rule. This suggests that in some cases, it might be beneficial to add more priority classes.
Theoretical analysis of the performance of SJF with predicted service times would be another interesting future research direction.
Lastly, it would be interesting to extend the analysis to multi-server queues (Some numerical experiments for multi-server queues is provided in Appendix \ref{app:multi}). When having multiple servers, we may be able to strike a better balance between system efficiency (average waiting time) and fairness among jobs of different sizes.\\

\begin{appendix}
\section{Proofs of the results in Section \ref{sec:main_proof}} \label{app:proof_prop}
We start by defining a few more notations. Let
\[\hat U_1^n(t)=\frac{K^n}{\sqrt{n}}U_1^n(nt) \mbox{ and } \hat Q_i^n(t)=\frac{K^n}{\sqrt{n}}Q_i^n(nt) \mbox{ for $i=1,2$}.\]

\subsection{Proof of of Proposition \ref{prop:class1}}
Let $W_k^n$ denote the waiting time of the $k$-th Class 1 arrival in queue in the $n$-th system.
Then based on Lindley's recursion for the single-server queue, we have for any $T>0$ \cite{BDP18},
\[\sup_{0\leq t \leq T}U_1^n(t) \leq \sup_{1\leq k\leq A_1^n(T)} (W_k^n+v_1^n(k)).\]
Thus, for any fixed constant $a>0$, we have
\[\begin{split}
\PP\left(\sup_{0\leq t\leq T} \hat U_1^n(t)>\epsilon\right)
\leq&\PP\left(\max_{1\leq k\leq A_1^n(nT)} \frac{K^n}{\sqrt{n}}(W_{k}^n+v_1^n(k))>\epsilon\right)\\
\leq&\underbrace{\PP\left(\max_{1\leq k\leq (1+a)nT} \frac{K^n}{\sqrt{n}}(W_{k}^n+v_1^n(k))>\epsilon\right)}_{\mbox{(I)}}+\underbrace{\PP(A_1^n(nT)>(1+a)nT)}_{\mbox{(II)}} .
\end{split}\]
We next analyze (I) and (II) one by one.

For (I), define $\Delta_k^n:=v_1^n(k-1)-\tau^n(k)$,
and $S_0^n=0$ and $S_k^n=S_{k-1}^n+\Delta_k^n$. Then,
$$W_{k+1}^n=\max\{W_k^n+\Delta_{k+1}^n,0\}=\max_{0\leq i\leq k+1}\{S_{k+1}^n-S_i^n\}
\overset{D}{=}\max_{0\leq i\leq k+1} S_i^n.$$
Let $M^n=\sup_{k\geq 1} S_{k}^n$. Then,
\[
\PP\left(\max_{1\leq k\leq N} W_{k}^n>a\right)
\leq\sum_{k=1}^{N} \PP(W_k^n>a) \leq N \PP(M^n>a).
\]
This further indicates that
\[\begin{split}
\PP\left(\max_{1\leq k\leq (1+a)nT} \frac{K^n}{\sqrt{n}}(W_{k}^n+v_1^n(k))>\epsilon\right)\leq&
\underbrace{\PP\left(\max_{1\leq k\leq (1+a)nT}\frac{K^n}{\sqrt{n}}v_1^n(k)>\epsilon/2\right)}_{\mbox{(i)}}\\
&+\underbrace{(1+a)nT\PP\left(\frac{K^n}{\sqrt{n}} M^n>\epsilon/2\right)}_{\mbox{(ii)}}.\end{split}\]
For (i), because $v_1^n(k)\leq K^n$ and $(K^n)^2/\sqrt{n}\rightarrow 0$ as $n\rightarrow\infty$ under Assumptions \ref{ass:service} and \ref{ass:threshold}, we have
(i) converges to 0 as $n\rightarrow\infty$.\\
For (ii), because
\[\E[\Delta_k^n]=\frac{1}{F(K^n)}\int_{0}^{K^n}xf(x)dx-\frac{1}{\lambda^nF(K^n)}=
\frac{1}{\lambda^nF(K^n)}(\rho_1^n-1)<0,\] 
$M^n$ is the all time maximum of a random walk with negative.
We define 
\[\psi_n(\theta):=\log \E[\exp(\theta \Delta_k^n)].\]
As $v_1^n(k-1)\leq K^n$, there exists $\epsilon>0$, such that $\psi_n(\theta)<\infty$ for any $\theta>-\epsilon$.
As $\psi_n^{\prime}(0)=\E[\Delta_k^n]<0$ and $\psi_n^{\prime\prime}(0)=\mbox{Var}(\Delta_k^n)>0$,
there exists $\theta_n>0$ such that $\psi_n(\theta_n)=0$. By Cram\'{e}r-Lundberg bound \cite{asmussen2008applied}, 
\[\PP(M^n>c) \leq \exp(-\theta_n c).\]
From Corollary 3 in \cite{GW94}, we also have
\[\theta_n=2\frac{1}{\mbox{Var}(\Delta_k^n)\lambda^nF(K^n)}(1-\rho_1^n)+o(\E[\Delta_k^n]).\]
Because $\mbox{Var}(\Delta_k^n)\rightarrow \sigma_a^2+\sigma_s^2$ and $\lambda^nF(K^n) \to 1$ as $n \to \infty$, 
\[\lim_{n \to \infty}\theta_n/(1-\rho_1^n)=2/(\sigma_a^2+\sigma_s^2).\]
This implies that $\theta_n\sim 1-\rho_1^n$.
Then,
\[\mbox{(ii)} \leq (1+a)nT\exp\left(-\theta_n\frac{\sqrt{n}\epsilon}{K^n2}\right)\rightarrow 0 \mbox{ as $n\rightarrow\infty$.}\]

For (II), let $A_0(t)$ denote the renewal process with interarrival times distributed as $\tau^{\infty}$.
Under Assumptions \ref{ass:arrive} and \ref{ass:arrive2}, by Theorem 1 in \cite{GW94b}, we have
\[\lim_{t\rightarrow\infty}\frac{1}{t}\log\E[\exp(\theta A_0(t))]=-\psi_{a}^{-1}(-\theta).\]
Then by Gartner-Ellis theorem, there exists $a>0$, such that
\begin{equation} \label{eq:ldp}
\frac{1}{t}\log \PP( A_0(t)>(1+a) t)=-I(a),
\mbox{ where } I(a)=\sup_{\theta}\{\theta a + \psi_{a}^{-1}(-\theta)\}>0.
\end{equation}
In particular, the probability
$\PP( A_0(t)>(1+a) t)$ decays exponentially fast in $t$.
This implies 
\[\PP(A_1^n(nT)>(1+a)nT) \leq \PP(A_0(nT)>(1+a)nT)\rightarrow 0 \mbox{ as $n\rightarrow\infty$.}\]

Above all, if $a>0$ is chosen according to \eqref{eq:ldp}, we have established that both (I) and (II) converge to zero as $\epsilon\rightarrow\infty$.
%
This concludes the proof of state-space collapse for $K^nU_1^n(nt)/\sqrt{n}$.

We next establish the state-space collapse for $K^nQ_1^n(nt)/\sqrt{n}$. It follows similar lines of arguments as
the proof of Theorem 4 in Section 3.2 of \cite{R84S}. However, there are several extra technical difficulties
we need to address along the way.
We denote $a_1^n(t)$ as the arrival time of the Class 1 customer that is in service at time $t$.
If there is no Class 1 customer in service, we set $a_1^n(t)=t$. We also define $\bar a_1^n(t)=a_1^n(nt)/n$.

As
\[A_1^n(t)-A_1^n(a_1^n(t))\leq Q_1^n(t)\leq A_1^n(t)-A_1^n(a_1^n(t))+1,\]
\begin{equation}\label{eq:bd1}
\frac{K^n}{\sqrt{n}}\left(A_1^n(nt)- A_1^n(n\bar a_1^n(t))\right)
\leq \hat Q_1^n(t)\leq \frac{K^n}{\sqrt{n}}\left(A_1^n(nt)-A_1^n(n\bar a_1^n(t))\right)+\frac{K^n}{\sqrt{n}}.
\end{equation}
We also note that as
\[U_1^n(a_1^n(t))\leq t-a_1^n(t)\leq U_1^n(a_1^n(t)) + v_1^n(A_1^n(a_1^n(t))),\]
\begin{equation}\label{eq:bd2}
\hat U_1^n(\bar a_1^n(t))\leq K^n\sqrt{n}(t-\bar a_1^n(t))\leq \hat U_1^n(\bar a_1^n(t)) + \frac{K^n}{\sqrt{n}}v_1^n(A_1^n(a_1^n(nt))).
\end{equation}
As $v_1^n(k)\leq K^n$,
\[\frac{K^n}{\sqrt{n}}v_1^n(A_1^n(a_1^n(nt)))\leq \frac{(K^{n})^2}{\sqrt{n}}\rightarrow 0 \mbox{ as $n\rightarrow\infty$},\]
where the convergence holds under Assumptions \ref{ass:service} and \ref{ass:threshold}.
Then, from \eqref{eq:bd2},
\begin{equation}\label{eq:time1}
\sup_{0\leq t\leq T}\left|K^n\sqrt{n}(t-\bar a_1^n(t))-\hat U_1^n(\bar a_1^n(t)) \right|\Rightarrow 0  \mbox{ as $n\rightarrow\infty$}.
\end{equation}

From \eqref{eq:bd1}, we have
\[\sup_{0\leq t\leq T}\hat Q_1^n(t)\leq \underbrace{\sup_{0\leq t\leq T}\frac{K^n}{\sqrt{n}}\left|A_1^n(nt)-A_1^n(n\bar a_1^n(t))\right|}_{\mbox{(A)}}+\frac{K^n}{\sqrt{n}}.\]
We first note under Assumptions \ref{ass:service} and \ref{ass:threshold}, $K^n/\sqrt{n}\rightarrow 0$ as $n\rightarrow0$.\\
%
Next, for (A), set $\delta_n=\frac{\epsilon}{4(1+a)}\frac{\sqrt{n}}{K^n}$, where $a$ is chosen according to \eqref{eq:ldp}. Then,
\[\begin{split}
&\PP\left(\sup_{0\leq t\leq T}\frac{K^n}{\sqrt{n}}\left|A_1^n(nt)-A_1^n(n\bar a_1^n(t))\right|>\epsilon \right)\\
=&\PP\left(\sup_{0\leq t\leq T}\frac{K^n}{\sqrt{n}}\left|A_1^n(nt)-A_1^n(n\bar a_1^n(t))\right|>\epsilon; \sup_{0 \le t \le T} |nt-n\bar a_1^n(t)| > \delta_n \right)\\
& + \PP\left(\sup_{0\leq t\leq T}\frac{K^n}{\sqrt{n}}\left|A_1^n(nt)- A_1^n(n\bar a_1^n(t))\right|>\epsilon; \sup_{0 \le t \le T} |nt-n\bar a_1^n(t)| \le \delta_n\right) \\
\leq & \underbrace{\PP\left(\sup_{0 \le t \le T} \frac{K^n}{\sqrt{n}} |nt-n\bar a_1^n(t)| > \frac{\epsilon}{4(1+a)} \right)}_{\mbox{(a)}}
+ \underbrace{\PP\left(\sup_{0\leq t\leq nT}\frac{K^n}{\sqrt{n}}\left|A_1^n(t+\delta_n)- A_1^n(t)\right|>\frac{\epsilon}{2}\right)}_{\mbox{(b)}}
\end{split}\]
For (a), from \eqref{eq:time1}, we have 
\[\frac{K_n}{\sqrt{n}}(nt-n\bar a_1^n(t))\Rightarrow \zeta(t) \mbox{ in $D[0,\infty)$ as $n\rightarrow\infty$}.\] 
Thus, (a) converges to zero as $n\rightarrow\infty$.
For (b), 
\begin{eqnarray*}
&&\PP\left(\sup_{0\leq t\leq nT}\frac{K^n}{\sqrt{n}}A_1^n(t+\delta_n)- A_1^n(t)>\frac{\epsilon}{2}\right)\\
&\le& \PP\left(\sup_{i=0, \dots, \lceil\frac{nT}{\delta_n}\rceil}A_1^n(i\delta_n+2\delta_n)- A_1^n(i\delta_n)>\frac{\epsilon}{2}\frac{\sqrt{n}}{K^n}\right)\\  
&\le& \left\lceil \frac{nT}{\delta_n}\right\rceil \PP\left(A_1^n(2\delta_n) > \frac{\epsilon}{2}\frac{\sqrt{n}}{K^n}\right)\\
&\leq& \left\lceil \frac{nT}{\delta_n}\right\rceil \PP\left(A_0(2\delta_n) > (1+a)2\delta_n\right)\\
&\leq& \left\lceil \frac{nT}{\delta_n}\right\rceil \exp(-\delta_n I(a)+\epsilon) \mbox{ from \eqref{eq:ldp} for $n$ large enough}\\
&\rightarrow& 0 \mbox{  as $n\rightarrow\infty$.}
\end{eqnarray*}
%
%
This completes the proof of the second part of Proposition \ref{prop:class1}.

\subsection{Proof of of Proposition \ref{prop:class2}}
Define 
\[\hat V_1^n(t)=\frac{1}{\sqrt{n}}\left(V_1^n(nt)+nt(1-\rho_1^n)\right).\]
Then,
\[\hat B^n(t):=\inf\{s\geq 0: \hat V_1^n(t+s/\sqrt{n}) - \hat V_1^n(t) + \hat U^n(t) \leq (1-\rho_1^n)s\}.\]
Let $\tau_1^n(k)$ be the interarrival time between the $k$-th and $(k-1)$-th arrival in $A_1^n$.
Under Assumptions \ref{ass:service} -- \ref{ass:threshold}, 
we have $\E[\tau_1^n(k)]\rightarrow 1$ and Var$(\tau_1^n(k))\rightarrow \sigma_a^2$ as $n\rightarrow\infty$.
In addition, $\sup_{n\geq1}\E[\tau_1^n(k)^2]\leq \E[(\tau_1^\infty)^3]<\infty$.
Then, 
\[\hat V_1^n \Rightarrow \mbox{BM}(0, \sigma_a^2+\sigma_s^2) \mbox{ in $D[0,\infty)$ as $n\rightarrow\infty$,}\]
where BM$(\mu,\sigma^2)$ denote a Brownian motion with drift coefficient $\mu$ and diffusion coefficient $\sigma$ \cite{R84S}.
The convergence implies that for any $T>0$,
\begin{equation}\label{eq:conv1}
\PP\left(\sup_{0\leq t\leq T}\sup_{0\leq s\leq n^{1/2-\delta}/K^n}|\hat V_1^n(t+s/\sqrt{n})-\hat V_1(t)|\leq \epsilon\right)\rightarrow 1 \mbox{ as $n\rightarrow\infty$},
\end{equation}
where $\delta$ is chosen according to \eqref{eq:rho1},
i.e., $(1-\rho_1^n)n^{1/2-\delta}/K^n\rightarrow\infty$ as $n\rightarrow\infty$.

We also note assuming
\[\sup_{0\leq t\leq T}\sup_{0\leq s\leq n^{1/2-\delta}/K^n}|\hat V_1^n(t+s/\sqrt{n})-\hat V_1(t)|\leq \epsilon,\] 
if $\sup_{0\leq t\leq T}\hat B^n(t)\leq n^{1/2-\delta}/K^n$, 
\[\sup_{0\leq t\leq T}|(1-\rho_1^n)\hat B^n(t) - \hat U^n(t)| \leq \epsilon;\]
and if $\hat B^n(t)> n^{1/2-\delta}/K^n$ for some $t>0$,
\[\hat U^n(t)>(1-\rho_1^n) n^{1/2-\delta}/K^n-\epsilon.\]
Because $\hat U^n\Rightarrow \mbox{RBM}(-\beta, 1+\sigma^2)$ as $n\rightarrow\infty$, then, by \eqref{eq:rho1},
\[\PP\left(\sup_{0\leq t\leq T}\hat U^n(t)>(1-\rho_1^n)n^{1/2-\delta}/K^n-\epsilon \right) \rightarrow 0 \mbox{ as $n\rightarrow\infty$.}\]
Thus, the convergence in \eqref{eq:conv1} implies
\[\PP\left(\sup_{0\leq t\leq T} |(1-\rho_1^n)\hat B^n(t) - \hat U^n(t)| \leq \epsilon \right)\rightarrow 1 \mbox{ as $n\rightarrow\infty$}.\] 
In addition, because $\rho_2^n/(1-\rho_1^n)\rightarrow 1$ as $n\rightarrow\infty$, we have
\[\rho_2^n \hat B^n(t) - \hat U^n(t) \Rightarrow \zeta(t) \mbox{ in $D[0,\infty)$ as $n\rightarrow\infty$}.\]
This concludes the proof of the first part of Proposition \ref{prop:class2}.

We now prove the second part. Let $a_2^n(t)$ denote the time of the oldest class 2 customer in the system at time $t$. We also define $\bar a_2^n(t)=a_2^n(nt)/n$.
Then,
\begin{equation}\label{eq:bd21}
A_2^n(t)-A_2^n(a_2^n(t))\leq Q_2^n(t)\leq A_2^n(t)-A_2^n(a_2^n(t))+1
\end{equation}
and
\begin{equation}\label{eq:bd22}
B^n(a_2^n(t), U_n(a_2^n(t)-))\leq t-a_2^n(t) \leq B^n(a_2^n(t), U_n(a_2^n(t))).
\end{equation}
From \eqref{eq:bd22}, we have
\[\rho_2^n \hat B^n(\bar a_2^n(t)-)\leq \rho_2^n\sqrt{n}(t-\bar a_2^n(t))\leq \rho_2^n\hat B^n(\bar a_2^n(t)).\]
From the first part of the proposition, $\rho_2^n\hat B_n(t)\Rightarrow \mbox{RBM}(-\beta,1+\sigma^2)$.
Therefore, 
\begin{equation}\label{eq:partII}
\sup_{0\leq t\leq T}\left|\rho_2^n\hat B^n(\bar a_2^n(t))-\rho_2^n\sqrt{n}(t-\bar a_2^n(t))\right|\Rightarrow 0 \mbox{ as $n\rightarrow\infty$.}
\end{equation}
Because $\rho_2^n\sqrt{n}\rightarrow\infty$ as $n\rightarrow\infty$, we also have
$\bar a_2^n \Rightarrow \eta$ in $D[0,\infty)$ as $n\rightarrow\infty$.

From \eqref{eq:bd21}, we have
\[\frac{K^n}{\sqrt{n}}\left(A_2^n(nt)- A_2^n(n\bar a_2^n(t))\right)
\leq \hat Q_2^n(t)\leq 
\underbrace{\frac{K^n}{\sqrt{n}}\left(A_2^n(nt)- A_2^n(n\bar a_2^n(t))\right)}_{\mbox{(I)}} + \frac{K^n}{\sqrt{n}}.\]
Under Assumptions \ref{ass:service} and \ref{ass:threshold}, $K^n/\sqrt{n}\rightarrow 0$ as $n\rightarrow\infty$.
For (I), we first note that 
\[\begin{split}
&\frac{K^n}{\sqrt{n}}\left(A_2^n(nt)- A_2^n(n\bar a_2^n(t))\right)\\
=&\underbrace{K^n\sqrt{\lambda_2^n}\frac{1}{\sqrt{\lambda_2^n n}}(A_2^n(nt)-\lambda_2^n nt)}_{\mbox{(i)}}-
\underbrace{K^n\sqrt{\lambda_2^n}\frac{1}{\sqrt{\lambda_2^n n}}(A_2^n(n\bar a_2^n(t))-\lambda_2^n n\bar a_2^n(t))}_{\mbox{(ii)}}
+\underbrace{\gamma^n\sqrt{n}\rho_2^n(t-\bar a_2^n(t))}_{\mbox{(iii)}}.
\end{split}\]
We next analyze $(A_2^n(nt)-\lambda_2^n nt)/\sqrt{\lambda_2^n n}$.
Let $\tau_2^n(k)$ denote the interarrival time between the $k$-th and the $(k-1)$-th arrival in $A_2^n$.
Then, $\E[\tau_2^n(k)]=1/\lambda_2^n$ 
and $\mbox{Var}(\tau_2^n(k))=(1-\bar F(K^n))/(\lambda_2^n)^2 + \sigma_a^2/(\lambda^n\lambda_2^n)$.
We also define $\tilde A_2^n$ as a renewal process 
in which the interarrival time between the $k$-th and $(k-1)$-th arrival is $\lambda_2^n\tau_2^n(k)$.
Then $\tilde A_2^n(\lambda_2^nt)=\tilde A_2^n(t)$. 
For  $\tilde A_2^n$, the interarrival time has mean $\E[ \lambda_2^n\tau_2^n(k)]=1$ and variance
\[\mbox{Var}(\lambda_2^n\tau_2^n(k))=(\lambda_2^n)^2\frac{1-\bar F(K^n)}{(\lambda_2^n)^2} + (\lambda_2^n)^2\frac{\sigma_a^2}{\lambda^n\lambda_2^n}\rightarrow 1 \mbox{ as $n\rightarrow\infty$}.\]
In addition, 
\[\sup_{n\geq 1}\E[(\lambda_2^n\tau_2^n(k))^3]\leq \sup_{n\geq 1}\left(\bar F(K^n)^2-6\bar F(K^n) + 6\right)\E[(\tau^\infty)^2]<\infty.\]
Because $\lambda_2^n n\rightarrow\infty$ as $n\rightarrow\infty$, 
\begin{equation}\label{eq:fcl}
\frac{A_2^n(nt)-\lambda_2^n nt}{\sqrt{\lambda_2^n n}}=\frac{\tilde A_2^n(\lambda_2^n nt)-\lambda_2^n nt}{\sqrt{\lambda_2^n n}}
\Rightarrow \mbox{BM}(0,1) \mbox{ in $D[0,\infty)$ as $n\rightarrow\infty$.}
\end{equation}

Now, for (i), under Assumption \ref{ass:service} and \ref{ass:threshold}, 
$K^n\sqrt{\lambda_2^n}\rightarrow 0$ as $n\rightarrow\infty$. Then, from \eqref{eq:fcl}, we have 
(i)$\Rightarrow \zeta$ in $D[0,\infty)$ as $n\rightarrow\infty$.
Similarly, for (ii), because $\bar a_2 \rightarrow \eta$ and $K^n\sqrt{\lambda_2^n}\rightarrow 0$, 
(ii)$\Rightarrow \zeta$ in $D[0,\infty)$ as $n\rightarrow\infty$.
Lastly, because $\gamma^n\rightarrow\gamma$ as $n\rightarrow\infty$,
from the first part of Proposition \ref{prop:class2} and \eqref{eq:partII}, we have
(iii)$\Rightarrow \gamma$RBM$(-\beta,\sigma_a^2+\sigma_s^2)$ in $D[0,\infty)$ as $n\rightarrow\infty$.

Putting (i), (ii), and (iii) together, we have (I)$\Rightarrow \gamma$RBM$(-\beta,\sigma_a^2+\sigma_s^2)$ in $D[0,\infty)$ as $n\rightarrow\infty$.
This further implies that
\[\hat Q_2^n(t)\Rightarrow \hat Q(t) \mbox{ in $D[0,\infty)$ as $n\rightarrow\infty$.}\]

%

\section{Proof of Theorem \ref{th:limit_inter}} \label{app:thm2}
From Theorem \ref{th:main}, we have for any $x\geq 0$,
\begin{equation}\label{eq:lim1}
\lim_{t\rightarrow\infty} \lim_{n\rightarrow\infty}\PP\left(\frac{K^n}{\sqrt{n}}Q^n(nt)\leq x\right)=1-\exp\left(-\frac{2\beta}{\gamma(1+\sigma_s^2)}x\right).
\end{equation}
For the first result in Theorem \ref{th:limit_inter}, we essentially want to show that we can change the limit on the right-hand side of \eqref{eq:lim1}.

For the $n$-th system, as $\rho^n<1$ and the scheduling policy is non-idling,
the stationary distribution of $Q^n(t)$ is well defined.
In particular, under the two-class preemptive priority rule, 
the steady-state queue length satisfies \cite{SM66}
\[\E[Q^n(\infty)]=\rho_1^n+\frac{\rho_2^n}{1-\rho_1^n}+\frac{\lambda_1^{n,2}\E[(v_1^{n})^2]}{2(1-\rho_1^n)}
+\frac{\lambda^n\lambda_2^n\E[v^2]}{2(1-\rho_1^n)(1-\rho^n)}.\]
As $\frac{\sqrt{n}}{K^n}(1-\rho_1^n)\rightarrow\infty$,
$\sqrt{n}(1-\rho^n)\rightarrow\beta$, and $K^n\mu_2^n\rightarrow\gamma$
as $n\rightarrow\infty$,
we have
\begin{equation}\label{eq:lim2}
\frac{K^n}{\sqrt{n}}\E[Q^n(\infty)]
\rightarrow\frac{\gamma(1+\sigma_s^2)}{2\beta}
\mbox{ as $n\rightarrow\infty$.}
\end{equation}
This establishes the second convergence result in Theorem \ref{th:limit_inter}.
 
The convergence in \eqref{eq:lim2} also implies $\frac{K^n}{\sqrt{n}}Q^n(\infty)$ is tight. 
Then, the interchange of limit result then follows from \cite{GZ06}.

\section{Proof of the results in Section \ref{sec:error}} \label{app:error}
\subsection{Proof of Proposition \ref{th:bounded}}
As $\varepsilon\in[-M,M]$ with probability one, we have the following bounds for $\hat \mu_2^n$ and $\hat \rho_1^n$:
\[ \frac{1}{\hat \mu_2^n}= \E[v|\hat v>K^n]=\E[v|\hat v>K^n, v> K^n-M]\geq K^n-M,\]
\[\begin{split}
\frac{1}{\hat \mu_2^n}= \E[v|\hat v>K^n]
&= \frac{\E[v1_{\hat v>K^n,K^n-M<v<K^n+M}]+\E[v1_{v>K^n+M}]}{\PP(\hat v>K^n,K^n-M<v<K^n+M) + \PP(v>K^n+M) }\\
&\leq \E[v|\hat v>K^n,K^n-M<v<K^n+M] + \E[v|v>K^n+M]\\
&\leq (K^n+M) + \frac{3}{2}(K^n+M) \mbox{ by Assumption \ref{ass:service}},
\end{split}\]
and
\[\begin{split}
1-\hat \rho_1^n = 1-\lambda^n\E[v1_{\hat v\leq K^n}]&=1-\E[v1_{\hat v\leq K^n}] + \frac{\beta}{\sqrt{n}}\E[v1_{\hat v\leq K^n}]\\ 
&=\E[v1_{\hat v> K^n}] + \frac{\beta}{\sqrt{n}}\E[v1_{\hat v\leq K^n}]\\
&\geq (\hat \mu_2^n)^{-1} \PP(\hat v>K^n) \geq  (\hat \mu_2^n)^{-1}\bar F(K^n+M) 
\end{split}\]
As $M\in(0,\infty)$ is fixed, $1/\hat \mu_2^n \sim K^n$ and for $0<\delta<1/2$ chosen according to Assumption \ref{ass:threshold}
\[
n^{1/2-\delta}\hat\mu_2^n(1-\hat \rho_1^n)\geq n^{1/2-\delta} \bar F(K^n+M) \rightarrow\infty \mbox{ as $n\rightarrow\infty$.}
\]
We have thus checked conditions C1 and C2. The rest of the process level convergence proof follows exactly the same lines of argument as the proof of Theorem \ref{th:main}. We thus omit it here.




\subsection{Proof of Proposition \ref{th:class}.}
Let $\mathcal{C}_i^n$ denote the event that the job is classified as Class $i$, $i=1,2$, in system $n$. 
We first note that 
\[\begin{split}
1-\hat \rho_1^n = 1-\lambda^n\E[v1_{\mathcal{C}_1^n}]&=1-\E[v1_{\mathcal{C}_1^n}] + \frac{\beta}{\sqrt{n}}\E[v1_{\mathcal{C}_1^n}]\\ 
&=\E[v1_{\mathcal{C}_2}] + \frac{\beta}{\sqrt{n}}\E[v1_{\mathcal{C}_1}]\\
&\geq (\hat \mu_2^n)^{-1} \PP(\mathcal{C}_2^n) = (\hat \mu_2^n)^{-1}\left(p_{12}^n F(K^n) + p_{22}^n\bar F(K^n)\right)
\end{split}\]
Thus, 
\[
n^{1/2-\delta}\hat \mu_2^n(1-\hat \rho_1^n) \geq n^{1/2-\delta} \left(p_{12}^n F(K^n) + p_{22}^n\bar F(K^n)\right) \rightarrow\infty \mbox{ as $n\rightarrow\infty$,}
\]
where the convergence follows because $n^{1/2-\delta}F(K^n)\rightarrow\infty$, $n^{1/2-\delta}\bar F(K^n)\rightarrow\infty$, and
$\lim_{n\rightarrow\infty} p_{12}^n+p_{22}^n>0$ under the assumption that $p_{21}^n\leq a$ for some $0\leq a<1$. We have thus verified
condition C2.

We next check condition C1, i.e., whether $\hat \mu_2^n$ converges to zero.
When $p_{12}^n\leq C \bar F(K^n)$ for some $C>0$ and $p_{21}^n\leq a$,
\[
\frac{1}{\hat \mu_2^n} = \frac{p_{12}^n\int_{0}^{K^n}xf(x)dx + p_{22}^n\int_{K^n}^{\infty}xf(x)dx}{p_{12}^n F(K^n) + p_{22}^n \bar F(K^n)} 
\leq \frac{C \bar F(K^n)+\int_{K^n}^{\infty}xf(x)dx}{(1-a)\bar F(K^n)} \sim K^n
\] 
and
\[
\frac{1}{\hat \mu_2^n} \geq \frac{(1-a)\int_{K^n}^{\infty}xf(x)dx}{C \bar F(K^n) + \bar F(K^n)} \sim K^n.
\]
Therefore, $1/\hat\mu_2^n\sim K^n$.The rest of the process level convergence proof follows from the proof of Theorem \ref{th:main}. 

On the other hand, when $p_{12}^n/(K^n\bar F(K^n))\rightarrow\infty$ as $n\rightarrow\infty$ and $p_{21}^n\leq a$,
\[
\frac{1}{\hat \mu_2^n} = \frac{p_{12}^n\int_{0}^{K^n}xf(x)dx + p_{22}^n\int_{K^n}^{\infty}xf(x)dx}{p_{12}^n F(K^n) + p_{22}^n \bar F(K^n)} 
\rightarrow 1 \mbox{ as $n\rightarrow\infty$}. 
\]
In this case, the Class 2 queue scales as $\sqrt{n}$.


\subsection{Proof of Proposition \ref{th:est}}
For any fixed $M>0$, we first note 
\[\begin{split}
\frac{1}{\hat \mu_2^n}=\frac{\int_{0}^{\infty}tf(t)\bar\Phi(K^n-t)dt}{\int_{0}^{\infty}f(t)\bar\Phi(K^n-t)dt}
&\geq \frac{M\int_{M}^{\infty}f(t)\bar\Phi(K^n-t)dt}{\int_{0}^{\infty}f(t)\bar\Phi(K^n-t)dt}\\
&=M\left(1- \frac{\int_{0}^{M}f(t)\bar\Phi(K^n-t)dt}{\int_{0}^{\infty}f(t)\bar\Phi(K^n-t)dt}\right).
\end{split}\]
We also note
\[\frac{\int_{0}^{M}f(t)\bar\Phi(K^n-t)dt}{\int_{0}^{\infty}f(t)\bar\Phi(K^n-t)dt}
\leq\frac{\bar\Phi(K^n-M)\int_{0}^{M}f(t)dt}{\bar\Phi(0)\int_{K^n}^{\infty}f(t)dt}
\leq\frac{\bar\Phi(K^n-M)}{\bar\Phi(0)\bar F(K^n)}\rightarrow 0,\]
as $n\rightarrow \infty$ under Assumption \ref{ass:error2}.
Thus, for any $M>0$, there exists $N>0$ such that for $n>N$, $1/\hat \mu_2^n>M/2$. As $M$ can be set arbitrarily large,
we have $1/\hat \mu_2^n\rightarrow \infty$ as $n\rightarrow\infty$.

We next analyze the case where for any fixed $a>0$, $\lim_{t\rightarrow\infty} \bar \Phi(at)/\bar F(t)=0$. 
Fix $a\in(0,1)$. Then,
\[\begin{split}
\frac{1}{\hat \mu_2^n}=\frac{\int_{0}^{\infty}tf(t)\bar\Phi(K^n-t)dt}{\int_{0}^{\infty}f(t)\bar\Phi(K^n-t)dt}
&\geq \frac{(1-a)K^n\int_{(1-a)K^n}^{\infty}f(t)\bar\Phi(K^n-t)dt}{\int_{0}^{\infty}f(t)\bar\Phi(K^n-t)dt}\\
&=(1-a)K^n\left(1- \frac{\int_{0}^{(1-a)K^n}f(t)\bar\Phi(K^n-t)dt}{\int_{0}^{\infty}f(t)\bar\Phi(K^n-t)dt}\right)\\
&\geq (1-a)K^n\left(1-\frac{\bar \Phi(aK^n)}{\bar\Phi(0)\bar F(K^n)}\right). 
\end{split}\]
As $\lim_{t\rightarrow\infty} \bar \Phi(at)/\bar F(t)=0$, for $n$ large enough, $1/\hat \mu_2^n>(1-a)K^n/2$.\\
On the other hand,
\[\begin{split}
\frac{1}{\hat \mu_2^n}=\frac{\int_{0}^{\infty}tf(t)\bar\Phi(K^n-t)dt}{\int_{0}^{\infty}f(t)\bar\Phi(K^n-t)dt}
&=\frac{\int_{0}^{K^n}tf(t)\bar\Phi(K^n-t)dt+\int_{K^n}^{\infty}tf(t)\bar\Phi(K^n-t)dt}{\int_{0}^{K^n}f(t)\bar\Phi(K^n-t)dt+\int_{K^n}^{\infty}f(t)\bar\Phi(K^n-t)dt}\\
&\leq \frac{K^n\int_{0}^{K^n}f(t)\bar\Phi(K^n-t)dt}{\int_{0}^{K^n}f(t)\bar\Phi(K^n-t)dt}
+\frac{\int_{K^n}^{\infty}tf(t)dt}{\bar \Phi(0)\int_{K^n}^{\infty}f(t)dt}\\
&=\left(1+\frac{1}{\bar\Phi(0)\gamma^n}\right)K^n
\leq \left(1+\frac{3}{2\bar\Phi(0)}\right)K^n \mbox{ by Assumption \ref{ass:service}}.
\end{split}\]
Thus, $1/\hat \mu_2^n\sim K^n$.

The rest of the proof follows exactly the same lines of argument as the proof of Theorem \ref{th:main}.

\subsection{Proof of Lemma \ref{lem:heavy}}
Note that
\[\begin{split}
\hat\mu_2^n&=\frac{\int_{0}^{\infty}f(t)\bar\Phi(K^n-t)dt}{\int_{0}^{\infty}tf(t)\bar\Phi(K^n-t)dt}\\
&=\frac{\int_{0}^{K^n/2}f(t)\bar\Phi(K^n-t)dt+\int_{K^n/2}^{\infty}f(t)\bar\Phi(K^n-t)dt}{\int_{0}^{K^n/2}tf(t)\bar\Phi(K^n-t)dt+\int_{K^n/2}^{\infty}tf(t)\bar\Phi(K^n-t)dt}\\
&\geq \frac{\bar\Phi(K^n)F(K^n/2)+\bar\Phi(K^n/2)\bar F(K^n/2)}{\bar\Phi(K^n/2) + \int_{K^n/2}^{\infty} tf(t)dt}\\
&\geq \frac{\bar\Phi(K^n)F(K^n/2)+\bar\Phi(K^n/2)\bar F(K^n/2)}{\bar\Phi(K^n/2) + K^n\bar F(K^n/2)}\\
&= \frac{\frac{\bar\Phi(K^n)}{\bar\Phi(K^n/2)}F(K^n/2)+\bar F(K^n/2)}{1 + \frac{K^n\bar F(K^n/2)}{\bar\Phi(K^n/2)}}\rightarrow \frac{1}{\xi(1/2)(1+2c)} \mbox{ as $n\rightarrow\infty$}.
\end{split}\]
As $(\xi(1/2)(1+2c))^{-1}>0$, this completes the proof of Lemma \ref{lem:heavy}.

\section{Exact analysis of different scheduling rules}
For $M/GI/1$ queues, we have closed-form expression for the steady-state average queue length
under five different scheduling rules \cite{SM66}: 1) two-class priority rule without preemption (NP);
2) two-class priority rule with preemption (P); 3) two-class priority rule with semi-preemption, where the priority of the job depends on its remaining size, i.e., if the remaining size of a job is smaller than $K^n$, then it is classified as a Class 1 job, otherwise, it is classified as a Class 2 job;
4) SJF; and 5) SRPT. We use $Q^n(\infty;i)$ to denote the stationary number of customers in the system under scheduling policy $i$, $i=1,2,\dots,5$. Assume $\mu=1$.
Under Two-Class NP,
\[\E[Q^n(\infty;1)]=\rho^n+\frac{\lambda^{n}\lambda_1^n\E[v^2]}{2(1-\rho_1^n)}
+\frac{\lambda^{n}\lambda_2^n\E[v^2]}{2(1-\rho_1^n)(1-\rho^n)}.\]
Under Two-Class P,
\[\E[Q^n(\infty;2)]=\rho_1^n+\frac{\rho_2^n}{1-\rho_1^n}+\frac{(\lambda_1^{n})^2\E[(v_1^{n})^2]}{2(1-\rho_1^n)}
+\frac{\lambda^n\lambda_2^n\E[v^2]}{2(1-\rho_1^n)(1-\rho^n)}.\]
Under Two-Class semi-preemptive, 
\[\E[Q^n(\infty;3)]=\rho_1^n + \frac{\rho_2^n}{1-\rho_1^n} - \frac{\lambda_2^n\rho_1^n K^n}{1-\rho_1^n}
+\frac{\lambda_1^n\lambda_2^n (K^{n})^2}{2(1-\rho_1^n)}
+\frac{(\lambda_1^{n})^2\E[(v_1^{n})^2]}{2(1-\rho_1^n)}
+\frac{\lambda^n\lambda_2^n\E[v^2]}{2(1-\rho_1^n)(1-\rho^n)}.\]
Under SJF,
\[\E[Q^n(\infty;4)]=\rho^n + \frac{(\lambda^{n})^2\E[v^{2}]}{2}\int_{0}^{\infty}\frac{1}{\left(1-\lambda^n\int_{0}^{x}ydF(y)\right)^2}dF(x).\]
Under SRPT,
\[\E[Q^n(\infty;5)]=\lambda^n\int_{0}^{\infty}\int_{0}^{x}\frac{1}{1-\lambda^n\int_{0}^{y}udF(u)}dydF(x)+\frac{(\lambda^{n})^{2}}{2}\int_{0}^{\infty}\frac{\int_{0}^{x}y^2dF(y)+\bar F(x)x^2}{\left(1-\lambda^n\int_{0}^{x}ydF(y)\right)^2}dF(x).\]

For the first three scheduling policies, $i=1,2,3$, under Assumption \ref{ass:threshold}, simple calculation reveals that
\[\frac{K^n}{\sqrt{n}}\E[Q^n(\infty;i)]\rightarrow\frac{\gamma(1+\sigma^2)}{2\beta}
\mbox{ as $n\rightarrow\infty$,}\]
or equivalently,
\[\frac{1}{\sqrt{n}\mu_2^n}\E[Q^n(\infty;i)]\rightarrow\frac{(1+\sigma^2)}{2\beta}
\mbox{ as $n\rightarrow\infty$.}\]
These limits agree with what we derived in Theorem \ref{th:limit_inter}, 
and suggests that preemption does not have an effect in the limit. In particular, 
the two-class priority rules achieve the same asymptotic performance under different preemption assumptions.

When comparing the two-class priority rule to SJF and SRPT, we note that
\[\E[Q^n(\infty;1)] \geq \E[Q^n(\infty;4)] \geq \E[Q^n(\infty;5)].\]
To quantify the actual differences between the three policies, we need to restrict ourselves to specific forms of service time distributions (see, e.g., \cite{LWZ11}). 


\section{Additional Numerical Experiments}
In this section, we provide some additional numerical experiments for the performance of the two-class priority rule.
In particular, we compare the two-class priority rule to another benchmark policy: foreground-background (FB) policy, which does not require any prior job-size information. We also study the performance of the two-class priority rule in multi-server queues.

\subsection{Comparison to the foreground-background policy}
Under FB, priority is given to the job that has received the least amount of service. If there are multiple such jobs, then they are served in a processor sharing manner \cite{NW07}. FB has several slightly different implementations. We choose this one as it has been shown to be optimal under certain conditions \cite{NW07}. Note that process-sharing is in general not applicable in service operations applications, and is thus not considered in rest of this paper.

FB is known to perform well for service time distributions with decreasing failure rates (heavy-tailed service time distributions). However, when the service time distribution
has an increasing failure rate, it can perform worse than FCFS.
In Table \ref{tab:FB}, we compare the performance of our two-class priority rule using predicted service time to the FB policy for $M/GI/1$ queues. We test two different service time distributions: one has a decreasing failure rate (Pareto$(m_{2.5},2.5)$) and the other has an increasing failure rate (Weibull$(\nu_{1.5},1.5)$). 
We observe that when the service time distribution has a decreasing failure rate, FB performs much better than FCFS but worse than our two-class priority rule. This observation is well-expected because the two-class priority rule utilizes more job-size information than FB. When the service time distribution has an increasing failure rate, FB actually performs worse than FCFS, while our two-class priority rule still performs very well, much better than FCFS. This observation suggests that the two-class priority rule is more robust to different service time distributions.
We also emphasize that when predicted service time information is available, the two-class priority rule is much easier to implement in service systems than FB. 

\begin{table}[htp]
\centering
\caption{Steady-state average queue length for an $M/GI/1$ queue under different scheduling policies.
For the two-class priority rules (P and NP), we use predicted service times where $\hat v = v+\epsilon$ and $\epsilon$'s are iid $N(0,0.3^2)$. We set the threshold $K(\rho,F)=\bar F((1-\rho)^{-(1-0.05)})$.
For the Pareto distribution, $\bar F(x)=(m_{2.5}/x)^{2.5}1\{x\geq m_{2.5}\}$. For the Weibull distribution, $\bar F(x)=\exp(-(x/\nu_{1.5})^{1.5})1\{x\geq 0\}$.} \label{tab:FB}
\begin{tabular}{c|cccc||cccc}
\hline
& \multicolumn{4}{|c||}{Pareto}   & \multicolumn{4}{|c}{Weibull} \\\hline
$\rho$       & 0.8   & 0.85 & 0.9  & 0.95   & 0.8   & 0.85 & 0.9  & 0.95   \\\hline
FCFS & 3.68 & 5.19 & 8.19 & 17.20    & 3.14 & 4.37 & 6.82 & 14.14      \\
Two-Class NP  & 2.87  & 3.72 & 5.18 & 8.64   & 2.42  & 3.14 & 4.46 & 8.01  \\
Two-Class P& 2.53 & 3.23 & 4.43 & 7.29 &  2.29  & 2.98 & 4.26 & 7.71       \\
FB & 3.22  & 4.27 & 6.12 & 10.43 &  4.64  & 6.77 & 11.20 & 25.30       \\\hline

\end{tabular}
\end{table}

\subsection{Multi-server queues} \label{app:multi}
As discussed in Section \ref{sec:lit}, even though we only prove the scaling results for single-server queues,
we expect the same results to hold for multi-server queues under the conventional heavy-traffic scaling. 
In this section, we analyze the performance of the two-class priority rule with nonpreemption in multi-server queues. 
We choose FCFS and SJF as two benchmark policies.

In Table \ref{tab:multi}, we compare the steady-state average queue length for $M/M/c$ queues under different scheduling policies.
We set $\mu=1$ and vary the traffic intensity, $\rho=\lambda/(c\mu)$, by varying the value of the arrival rate $\lambda$. 
We observe that for large values of $c$, when $\rho$ is small, there is almost no performance improvement from ``smart" scheduling. 
The reason is that these systems require almost no waiting.
In particular, the system operates like an $M/M/\infty$ queue, in which all customers enter service immediately upon arrival. Thus, no scheduling decision needs to be made. However, when $\rho$ is large enough, we still gain significant performance improvement from job-size-based scheduling rules. In addition, our two-class priority rule achieves almost the same performance as SJF, but is much easier to administer in practice.


\begin{table}[htp]
\centering
\caption{Steady-state average queue length for $M/M/c$ queues under different scheduling policies ($\mu=1$,$\lambda=\rho c$, $K(\rho,F)=\bar F((1-\rho)^{(1-0.05)})$).} \label{tab:multi}
\begin{tabular}{c|ccc||ccc}
\hline
& \multicolumn{3}{c||}{$c=10$}  & \multicolumn{3}{c}{$c=100$}\\\hline

$\rho$    & 0.80          & 0.90           & 0.99      & 0.80          & 0.90           & 0.99 \\ \hline
FCFS     & 9.64         & 15.02         & 105.31        & 80.0 & 92.0 & 186.4        \\
Two-Class NP  & 8.96 & 11.87        & 37.74          & 80.0 & 90.9  & 124.5      \\
SJF       & 8.77       & 11.14          &  26.72         & 80.0  & 90.7  & 114.5  \\\hline
\end{tabular}
\end{table}

\end{appendix}

\bibliographystyle{plain}
\bibliography{myrefsBounding}
\end{document}